\documentclass[11pt]{article}
\usepackage{graphicx}
\usepackage{amsmath}
\usepackage{amsthm}
\usepackage{amssymb}
\usepackage[letterpaper, top=1in, bottom=1in, left=1.25in, right=1.25in]{geometry}
\newtheorem{theorem}{Theorem}
\newtheorem{lemma}[theorem]{Lemma}
\newtheorem{corollary}[theorem]{Corollary}

\newcommand{\mc}{\mathcal}

\newcommand{\spts}{\textsc{S}}
\newcommand{\ps}{\textsc{P}}
\newcommand{\es}{\textsc{E}}
\newcommand{\conv}{\mathrm{conv}}

\title{Visibility Sets, Empty Pentagons, and\\ Sufficient Conditions for Hexagons}
\author{Carlos M. Nicol\'as}
\date{}

\begin{document}	
\maketitle
\thispagestyle{empty}

\begin{abstract}
	Given a set $S$ of $n$ points in general position in the plane, we define the $C$-set $C(p,q)$ of an oriented edge $(p,q)$ as the set of extreme points of the visibility set $V(p,q)$ facing $(p,q)$. Since every point of $V(p,q)$ that is not in $C(p,q)$ gives rise to an empty convex pentagon having $(p,q)$ as a boundary edge, it is natural to study the properties of $C$-sets in order to bound the number of empty convex pentagons of $S$.
	
	We prove tight bounds on the intersection size and union complexity of $C$-sets sharing an endpoint. We also consider the problem of bounding the number of empty pentagons and hexagons in terms of the number of interior and boundary points of $S$. While the hexagon bound is weak, it offers a potential step toward a superlinear lower bound. On the other hand, the two-parameter bound for pentagons, combined with the structural results on $C$-sets, yields a simple proof of an $\Omega(n^{3/2})$ bound for the number of pentagons. While this is a significant improvement over previously published bounds, a recent preprint \cite{Astudillo2026} attains a better bound using completely different methods. Thus, the primary contribution of this paper lies in the structural results about $C$-sets and the two-parameter bounds for empty pentagons and hexagons.
\end{abstract}

\section{Introduction}

Let $S$ be a set of points in the plane in general position. A convex polygon with vertices in $S$ that does not contain any vertices of $S$ in its interior is called an empty convex polygon of $S$. A classical problem in discrete geometry, formulated by Erd\H{o}s in 1978~\cite{Erdos1978}, asks for bounds on the number of empty convex polygons determined by $S$. To be precise, let $g_k(S)$ be the number of empty convex polygons of $S$ with $k$ vertices, and let $f_k(n)$ be the minimum value of $g_k(S)$ among all point configurations $S$ of size $n$ in general position. The question posed by Erd\H{o}s asked whether $f_k(n)>0$ for sufficiently large values of $n$. This is obviously true for $k=3$ and $k=4$. 

For $k=5$, Harborth~\cite{Harborth1978} proved that $f_5(n)>0$ for $n\geq 10$. For $k=6$, the existence of empty hexagons was settled (almost 30 years later) in the affirmative by Nicol\'as~\cite{Nicolas2007} and independently by Gerken~\cite{Gerken2008}. On the other hand, Horton~\cite{Horton1983} showed the existence of arbitrarily large point configurations containing no empty heptagons, i.e., $f_7(n)=0$ for all $n$, which implies $f_k(n)=0$ for all $k\geq 7$ and all $n$. 

Given that $f_k(n)>0$ for $3\leq k\leq 6$ and sufficiently large $n$, it is natural to ask for the asymptotic rate of growth of $f_k$ for these values of $k$. B\'ar\'any and Valtr showed quadratic upper bounds for all $f_3$, $f_4$, $f_5$, and $f_6$~\cite{Barany2004}, and this is conjectured to be the actual rate of growth for all these functions. In the case of $f_3$ and $f_4$, quadratic lower bounds are known, so the problem is essentially settled for these two functions, though there is still a gap in the coefficients of the quadratic terms in the lower and upper bounds (see~\cite{Aichholzer2014} for the current best lower bounds). For $f_6(n)$, the best known lower bound is linear (see~\cite{Valtr2012}). That was also the case for $f_5$ until the year 2020, when Aichholzer, Balko, et al.~\cite{Aichholzer2020} obtained the first superlinear lower bound, $f_5(n)=\Omega(n\log^{4/5}(n))$, though it should be mentioned that parts of their proof are computer-assisted.

Recently, a preprint by Astudillo-Marb\'{a}n and Sol\'{e}-Pi~\cite{Astudillo2026} was posted on arXiv with an $\Omega(n^{20/11})$ lower bound. In our paper, we present a completely independent proof establishing an $\Omega(n^{3/2})$ lower bound. While the asymptotic growth rate of our bound is superseded by the recent preprint, we believe that some of the methods and results in this paper have intrinsic interest that extends beyond the problem of counting empty pentagons. For example, we study the set of edge-facing extreme points of the visibility sets of oriented edges, which we call $C$-sets, and establish tight bounds for the size of their intersection when the edges have a common endpoint. Moreover, also for edges with common endpoints, we show that the union of their $C$-sets can be partitioned into a linear number of angle-disjoint intervals (as seen from the common endpoint). Additionally, one of the main results in this paper (Theorem~\ref{pentagons_in_convex_set_Theorem}) gives bounds on the number of empty pentagons and hexagons in terms of the number $m$ of points on the boundary of $S$ and the number $n_0$ of points in the interior of $S$. For hexagons, the bound is quadratic in $m$ provided that the ratio $n_0/m$ is less than a certain constant. This could potentially be of use to prove a non-linear lower bound on the number of empty hexagons in the general case. For pentagons, the two-parameter bound, in conjunction with the results about $C$-sets, yields a simple proof of an $\Omega(n^{3/2})$ lower bound for the number of empty pentagons.

\begin{figure}[tpb!]
	\centering\includegraphics[scale=.7]{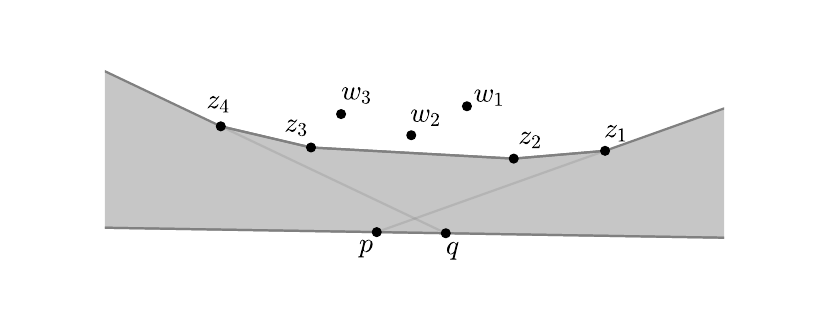}
	\caption{All points $z_i$ and $w_i$ belong to $V(p,q)$, but only the points $z_i$ belong to $C(p,q)$. The shaded region is empty.}
	\label{Figure11}
\end{figure}

\section{Visibility Sets and an Extension of Harborth's Theorem}
Throughout this paper, $\spts$ denotes a set of points in general position in $\mathbb{R}^2$, and $n$ denotes its cardinality. 
A point $w$ is a convex combination of the points $p$ and $q$ if there exist real numbers $\alpha\geq 0$ and $\beta\geq 0$ such that $w=\alpha p+\beta q$ and $\alpha+\beta=1$. The set of convex combinations of two points $p$ and $q$ is called the segment with endpoints $p$ and $q$. 
A set $\mc{R}\subseteq \mathbb{R}^2$ is convex if it contains the segment with endpoints $p$ and $q$, for all $p$ and $q$ in $\mc{R}$.
The convex hull of $\mc{R}\subseteq \mathbb{R}^2$, denoted $\conv(\mc{R})$, is the smallest convex set in $\mathbb{R}^2$ containing $\mc{R}$. We write $\conv(z_1,\ldots,z_k)$ instead of $\conv(\{z_1,\ldots,z_k\})$ for brevity.

A set $T\subseteq \spts$ is in convex position if all its elements lie on the boundary of its convex hull. If $T\subseteq \spts$ is in convex position, then $\conv(T)$ is called a convex polygon of $\spts$ with $|T|$ vertices. If $T\subseteq \spts$ is in convex position and $\conv(T)\cap \spts=T$, then $\conv(T)$ is called an empty convex polygon of $\spts$. In other words, a convex polygon is empty if it does not contain any points of $\spts$ in its interior. Naturally, an empty convex polygon with five vertices is called an empty convex pentagon. For the sake of brevity we occasionally omit the qualifiers `empty' and `convex' if it is clear from the context. Note that empty convex polygons with $k$ vertices are sometimes called $k$-holes in the literature (e.g., \cite{Aichholzer2014}). 
The set of all empty convex pentagons in $\spts$ will be denoted by $\ps$ throughout this paper. Our goal is to show that there exist constants $K>0$ and $N_0$ such that $|\ps|> Kn^{3/2}$ for all $n> N_0$. 

For any $\{p,q\}\subseteq \spts$ with $p\neq q$, the oriented edge from $p$ to $q$ is denoted $(p,q)$. All edges in this paper are oriented edges, even when we omit the word `oriented'. We refer to the side to the left of $(p,q)$, i.e., the side in which the ray from $p$ through $q$ falls if rotated in the counterclockwise direction, as the positive side of $(p,q)$. The other side is the negative side. 
The set of all oriented edges of $\spts$ is denoted $\es$, i.e.,  $\es=\{(p,q): \{p,q\}\subseteq \spts \mbox{ and }p\neq q\}$. The oriented line through $p$ and $q$ is denoted as $\mc{L}(p,q)$. The open half-planes to the left and right of $\mc{L}(p,q)$ are denoted by $\mc{H}^+(p,q)$ and $\mc{H}^-(p,q)$, respectively. We say that a line separates two sets if one set lies to the left of the line and the other one to the right. The set of pentagons containing $\conv(p,q)$ on their boundary and which are contained in $\mc{L}(p,q)\cup\mc{H}^+(p,q)$ is denoted $\ps(p,q)$, i.e., $\ps(p,q)=\{Q\in\ps: Q=\conv(p,q,z_1,z_2,z_3) \mbox{ and }\{z_1,z_2,z_3\}\subset\mc{H}^+(p,q)\}$.

For any $(p,q)\in\es$, the visibility set $V(p,q)$ of the edge $(p,q)$ is the set of points $w$ such that  $w\in \spts\cap\mc{H}^+(p,q)$ and $\conv(p,q,w)$ is empty. If $w\in V(p,q)$, we say that $w$ is visible from  $(p,q)$ and that $(p,q)$ is visible from $w$.

Let $C(p,q)$ denote the set of points $z\in V(p,q)$ such that $\conv(p,q,z)$ and  $\conv(\spts\cap\mc{H}^+(p,q))$ have disjoint interiors. 
Equivalently, $z\in C(p,q)$ if and only if $z\in V(p,q)$ and there exists a line through $z$ that separates $\spts\cap\mc{H}^+(p,q)-\{z\}$ from $\{p,q\}$. Note that $C(p,q)$ is in convex position since it lies on the boundary of $\conv(\spts\cap\mc{H}^+(p,q))$.
Note that if $C(p,q)=\{z_1, z_2, \ldots, z_k\}$, where $z_1, \ldots, z_k$ are listed in counterclockwise order as seen from $p$ (or from $q$), then $(\mc{H}^+(z_1,p)\cup\mc{H}^+(z_1,z_{2})\cup\cdots\cup\mc{H}^+(z_{k-1},z_{k})\cup\mc{H}^+(q,z_k))\cap\mc{H}^+(p,q)$ contains no points of $\spts$. We refer to this region as the empty region  between $\mc{L}(p,q)$ and the extended polygonal curve associated with $C(p,q)$.
See Figure~\ref{Figure11} for an example.

\begin{figure}[tpb!]
	\centering\includegraphics[scale=1]{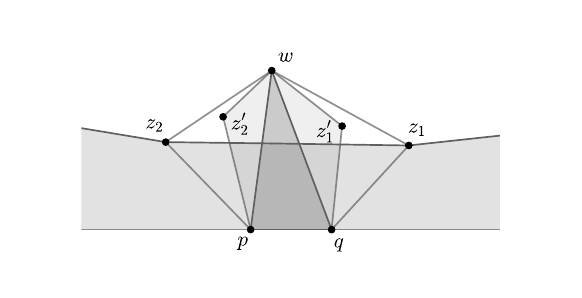}
	\caption{Point $w \in V(p,q) - C(p,q)$, while points $z_i \in C(p,q)$ for $i=1,2$. $\conv(p,q,z_1',w,z_2')$ is an empty convex pentagon. }
	\label{Figure6}
\end{figure}

\begin{theorem} For any distinct points $p$ and $q$ of $\spts$, $|V(p,q)|\leq |C(p,q)|+|\ps(p,q)|$.
\label{vcp_theorem}
\end{theorem}

\begin{proof} Let $w\in V(p,q)-C(p,q)$. Then, the triangle $\conv(p,q,w)$ is empty and intersects the interior of $\conv(\spts\cap\mc{H}^+(p,q))$. Therefore, $\conv(p,q,w)$ must intersect the boundary of $\conv(\spts\cap\mc{H}^+(p,q))$. Because $\conv(p,q,w)$ is empty, it cannot contain any vertices of $C(p,q)$; therefore, it must intersect the interior of a single edge $(z_1,z_2)$ on this boundary. Note that $z_1$ and $z_2$ are visible from $(p,q)$, so they belong to $C(p,q)$. Let $(z_1,z_2)$ be oriented so that $\spts\cap\mc{H}^+(p,q)-\{z_1,z_2\}$ is on the negative side of $(z_1,z_2)$. Under this orientation, $z_1$ belongs to $\mc{H}^+(w,q)$ and $z_2$ belongs to $\mc{H}^+(p,w)$. See Figure~\ref{Figure6}.

The empty triangle $\conv(p,q,w)$ can be extended to an empty pentagon $\conv(p,q,z_1', w, z_2')$ as follows: if $\conv(w,q,z_1)$ is empty, set $z_1'=z_1$; if not, choose $z_1'$ in the interior of $\conv(w,q,z_1)$ such that $\conv(w,q,z_1')$ is empty. $z_2'$ can be selected in a similar way. Since $(z_1,z_2)$ intersects $(p,w)$ and $(q,w)$, all five vertices $p,q,z_1',w,z_2'$  lie on the boundary of $\conv(p,q,z_1',w,z_2')$, guaranteeing its convexity. Furthermore, because it is the union of the empty triangles $\conv(p,q,w)$, $\conv(w,q,z_1')$, and $\conv(p,w,z_2')$, the resulting polygon $\conv(p,q,z_1',w,z_2')$ is an empty convex pentagon, i.e., it belongs to $\ps(p,q)$.

Finally, if both $w$ and $w'$ belong to $V(p,q)-C(p,q)$, with $w'\neq w$, then the empty pentagons obtained from $w$ and  $w'$ with this procedure must be distinct, since they differ on the vertices opposite to $(p,q)$. Therefore, $|\ps(p,q)|\geq |V(p,q)-C(p,q)|=|V(p,q)|-|C(p,q)|$.
\end{proof}

\begin{corollary} For any set of edges $D\subseteq \es$, $\sum_{(p,q)\in D}|V(p,q)|\leq 5|\ps|+\sum_{(p,q)\in D}|C(p,q)|$.
\label{vcp_corollary}
\end{corollary}

\begin{proof} By Theorem~\ref{vcp_theorem}, 
$\sum_{(p,q)\in D}|V(p,q)|\leq\sum_{(p,q)\in D}|\ps(p,q)|+\sum_{(p,q)\in D}|C(p,q)|\leq 5|\ps|+\sum_{(p,q)\in D}|C(p,q)|$, since any pentagon is counted at most five times over the distinct edges $(p,q)\in D$.
\end{proof}

Next, we establish a strengthening of Harborth's Theorem on the existence of empty pentagons. Specifically, we show that the pentagon either contains the top point of the set or contains an edge visible from the top point.

\begin{figure}[tpb!]
\centering\includegraphics[scale=1.25]{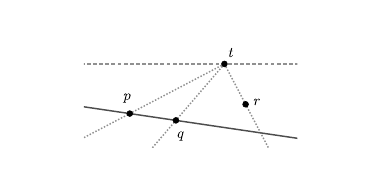}
\caption{$\angle pqr$ points away from $t$.}
\label{points_away}
\end{figure}

We refer to the point with the largest $y$-coordinate among the points in a set as the top point of the set. 
We say that $r$ lies above
 $\mc{L}(p,q)$ if $r\in\mc{H}^+(p,q)$, and $r$ lies below $\mc{L}(p,q)$ if $r\in\mc{H}^-(p,q)$. 
Also, if $p$, $q$, $r$ appear in that order when listed in counterclockwise order about a point $t$ and all three lie on the same side of a line through $t$, we say that $\angle pqr$ {points toward} $t$ if $r$ lies below $\mc{L}(p,q)$ and {points away from} $t$ if $r$ lies above $\mc{L}(p,q)$. See Figure~\ref{points_away} for an example.
Note that in the theorem below, when we refer to an edge of a pentagon we mean either a boundary edge or an internal diagonal of the pentagon. Thus, a pentagon has $20=2\binom{5}{2}$ oriented edges.

\begin{theorem} Let $A$ be a set of 10 points in general position and let $t$ be the top point of $A$. Then there exists an empty convex pentagon in $A$ that either contains $t$ or contains an edge that is visible from $t$. 
\label{GH_theorem}
\end{theorem}

\begin{figure}[tpb!]
	\centering $\includegraphics[scale=.7]{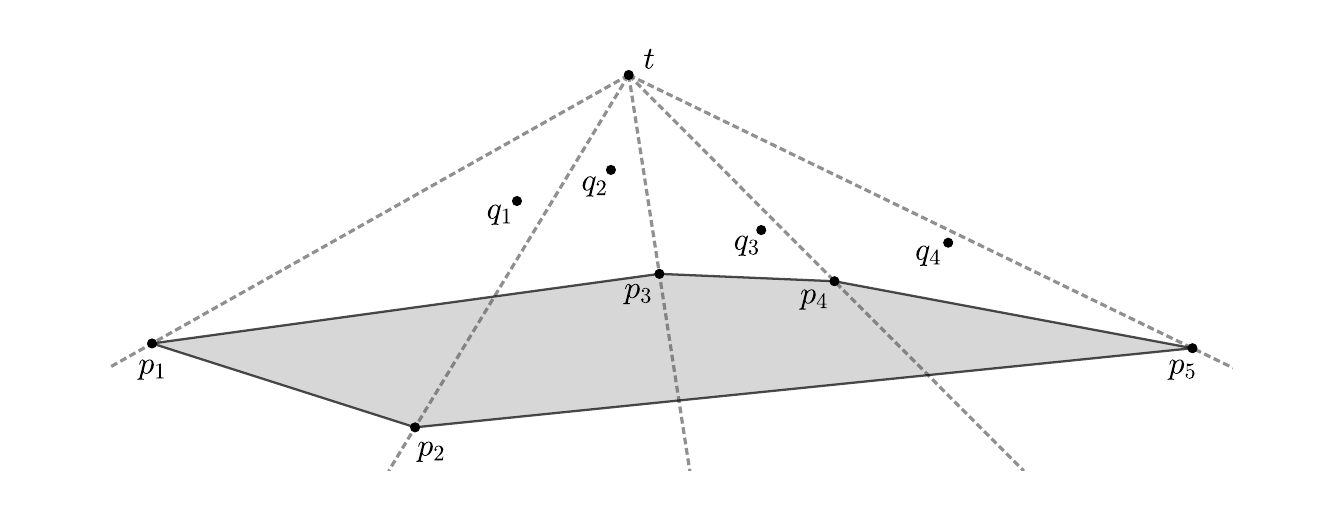}$
	\caption{$\conv(p_4, q_3, q_2, q_1, p_3)$ is an empty pentagon with edges visible from $t$.}
	\label{ten_points_figure}
\end{figure}

\begin{proof} 
	By Harborth's Theorem, $A$ contains an empty pentagon $Q$. If $t\in Q$, the theorem holds. Suppose instead that $t\not\in Q$. Let $p_1,\ldots, p_5$ be the vertices of $Q$ listed in counterclockwise order as seen from $t$. There are two paths from $p_1$ to $p_5$ along the boundary of $Q$. We call the path closer to $t$ the {upper boundary} of $Q$, and the path farther away the {lower boundary} of $Q$. Let $q_1,\ldots, q_4$ be the remaining four points in $A-\{t,p_1,\ldots,p_5\}$, listed in counterclockwise order as seen from $t$. 
	
	{Case 1: An edge of $Q$ is visible from $t$.} 
	If the triangle $\conv(t, p_i, p_{i+1})$ is empty for some $1\leq i\leq 4$, then the edge $(p_i, p_{i+1})$ of $Q$ is visible from $t$, satisfying the theorem. 
	
	{Case 2: No edge of $Q$ is visible from $t$.} 
	If Case 1 does not hold, then necessarily $q_i\in\conv(t, p_i, p_{i+1})$ for each $1\leq i\leq 4$. Note that this forces $\angle p_iq_ip_{i+1}$ to point toward $t$ for $1\leq i\leq 4$. 
	
	{Subcase 2.1: An angle $\angle q_i p_{i+1} q_{i+1}$ points toward $t$.} 
	If $\angle q_1p_2q_2$ points toward $t$ (which is possible when $p_2$ lies on the upper boundary of $Q$), then $\conv(p_3,q_2,p_2,q_1,p_1)$ is an empty pentagon with four edges visible from $t$. Similarly, if either $\angle q_2p_3q_3$ or $\angle q_3p_4q_4$ points toward $t$, we obtain an empty pentagon with edges visible from $t$. Therefore, the theorem holds in this case.
	
	{Subcase 2.2: } 
	Assume that $\angle q_ip_{i+1}q_{i+1}$ points away from $t$ for $1\leq i\leq 3$. Consider the line $\mc{L}(q_2,q_3)$. If $p_2$ lies above this line, then $\angle p_2q_2q_3$ points away from $t$ and $\conv(t,q_1,p_2,q_2,q_3)$ is an empty pentagon (since $\angle q_2p_3q_3$ points away from $t$) containing $t$. Similarly, if $p_4$ lies above $\mc{L}(q_2,q_3)$, then $\angle q_2q_3p_4$ points away from $t$ and $\conv(t,q_2,q_3,p_4,q_4)$ is an empty pentagon containing $t$. 
	
	Therefore, we may assume that both $p_2$ and $p_4$ lie below $\mc{L}(q_2,q_3)$. We now consider the position of $p_3$ on $Q$.
	
	{Subcase 2.2.1:} 
		Assume $p_3$ lies on the lower boundary of $Q$. Then $\angle p_2p_3p_4$ points away from $t$, and $\conv(p_2,p_3,p_4,q_3,q_2)$ is an empty pentagon with multiple edges visible from $t$.
		
	{Subcase 2.2.2:}	
		Assume $p_3$ lies on the upper boundary of $Q$. We examine the angles formed by the $q_i$ points. If both $\angle q_1q_2q_3$ and $\angle q_2q_3q_4$ point away from $t$, then $\conv(t,q_1,q_2,q_3,q_4)$ is an empty pentagon containing $t$.

		If not, at least one points toward $t$. Without loss of generality, assume that $\angle q_1q_2q_3$ points toward $t$. This implies $p_3$ lies below $\mc{L}(q_1,q_2)$, since $\angle q_2p_3q_3$ points away from $t$. We now consider the position of $p_2$ on $Q$:
		
	{Subcase 2.2.2.1:} 	
		Suppose $p_2$ lies on the upper boundary of $Q$ (like $p_3$): Consider the line $\mc{L}(p_2,p_3)$. If $q_1$ lies below this line, $\conv(p_4,p_3,p_2,q_1,p_1)$ is an empty pentagon with multiple edges visible from $t$. If $q_3$ lies below $\mc{L}(p_2,p_3)$, $\conv(p_4,q_3,p_3,p_2,p_1)$ is an empty pentagon with multiple edges visible from $t$. If both $q_1$ and $q_3$ lie above $\mc{L}(p_2,p_3)$, $\conv(q_1,p_2,p_3,q_3,q_2)$ is an empty pentagon (since $\angle q_1q_2q_3$ points toward $t$) with six edges visible from $t$. 
			
	{Subcase 2.2.2.2:} 
		Suppose $p_2$ lies on the lower boundary of $Q$: If $q_1$ lies below $\mc{L}(p_3,p_4)$, then $\conv(p_4,p_3,q_1,p_1,p_2)$ is an empty pentagon with two edges visible from $t$, since $p_3$ lies below $\mc{L}(q_1,q_2)$. If $q_1$ lies above $\mc{L}(p_3,p_4)$, then $\conv(p_4, q_3, q_2, q_1, p_3)$ is an empty convex pentagon (since $p_4$ lies below $\mc{L}(q_2,q_3)$ and $\angle q_1q_2q_3$ points toward $t$) with multiple edges visible from $t$ (this is the case depicted in Figure~\ref{ten_points_figure}).
\end{proof}

\begin{figure}[tpb!]
	\centering\includegraphics[scale=.4]{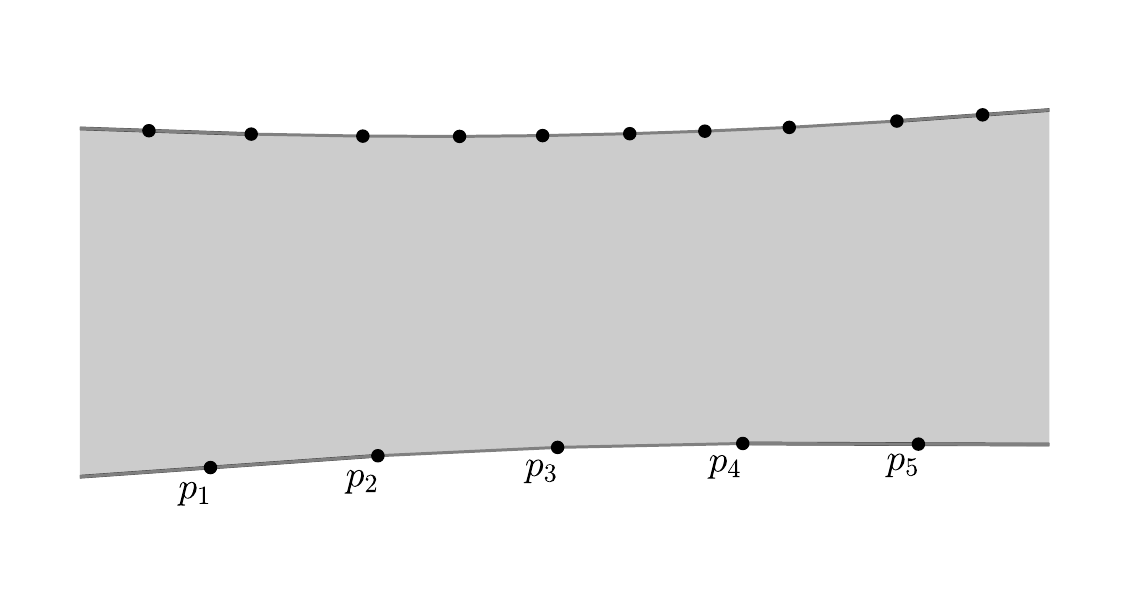}
	\caption{$C(p_i,p_{i+1})=C(p_j,p_{j+1})$ for all $i$ and $j$. The shaded area is empty.}
	\label{commonCset}
\end{figure}

We now we apply Theorem~\ref{GH_theorem} to a quadratic number of 10-point subsets. We treat each point in $\spts$ as a top point and angularly sweep through the points below it, forming wedges of 9 consecutive points. By Theorem~\ref{GH_theorem}, each of these wedges contributes either a pentagon containing its top point or a pentagon with an edge visible from it. Theorem~\ref{bound_on_edges_theorem} aggregates these local contributions.
	
\begin{theorem} Suppose $n\geq 10$. Then, there exists a set of edges $D\subseteq \es$ such that
	\begin{itemize}
		\item $|\ps|\geq \frac{1}{20}|D|$
		\item $6|\ps|+\sum_{(p,q)\in D}|C(p,q)|\geq \frac{1}{16}(n^2-17n)$
	\end{itemize}
\label{bound_on_edges_theorem}
\end{theorem}
	
\begin{proof} Let $v_1, \ldots, v_n$ be the $n$ points of $\spts$ sorted in increasing order by their $y$-coordinates. For $1\leq j<i\leq n$, define $w_j^i$ to be the $j$-th element of $\{v_1,\ldots,v_{i-1}\}$ when sorted in counterclockwise order about $v_i$. Let $I=\{(i,j): 10\leq i\leq n\mbox{ and }1\leq j\leq i-9\}$. For $(i,j)\in I$, let $\spts^i_j$ be the 9-element set $\{w^i_{j},w^i_{j+1},\ldots,w^i_{j+8}\}$. 
		
Applying Theorem~\ref{GH_theorem} to $\spts^i_j\cup\{v_i\}$, we obtain a pentagon $Q^i_j$ that is empty relative to $\spts^i_j\cup\{v_i\}$, and that either contains $v_i$ or contains an edge visible from $v_i$, relative to $\spts^i_j\cup\{v_i\}$. 
Now we show that these properties hold with respect to the entire set $\spts$. Since $v_i$ has the largest $y$-coordinate among the points in $\spts^i_j\cup\{v_i\}$, all of $\conv(\spts^i_j\cup\{v_i\})$ must lie in the closed half-plane below the horizontal line through $v_i$. Also, because $\spts^i_j$ is formed by 9 consecutive points in the angular sorting around $v_i$, there are no other points of $\{v_1, \dots, v_{i-1}\}$ in this wedge. 
Therefore,  $\spts\cap\conv(\spts^i_j\cup\{v_i\})=\spts^i_j\cup\{v_i\}$, so the properties of $Q^i_j$ (being empty and possessing a visible edge) hold with respect to $\spts$ as well. In particular, $Q^i_j\in \ps$. 
		
Now let $I_1=\{(i,j)\in I: v_i\in Q^i_j\}$ and $I_2=I-I_1$. For $(i,j)\in I_2$, define $p^i_j$ and $q^i_j$ to be the vertices of $Q^i_j$ such that the edge $(p^i_j,q^i_j)$ is visible from $v_i$. Finally, let $D=\{(p^i_j,q^i_j): (i,j)\in I_2\}$. 
Clearly, $|\ps|\geq \frac{1}{20}|D|$ because the edges in $D$ come from pentagons in $\ps$ and there are $2\binom{5}{2}=20$ oriented edges in any pentagon. 
		
To obtain a bound relating the size of the visibility sets of the edges in $D$ with the size of $|I_2|$, let us consider the function $f:I_2\rightarrow D$, defined by $f(i,j)=(p^i_j,q^i_j)$. For a fixed top point $v_i$, if two different indices $j_1$ and $j_2$ map to $(p,q)$, then both $\spts^i_{j_1}$ and $\spts^i_{j_2}$ contain $\{p,q\}$. Since $\spts^i_{j_1}$ and $\spts^i_{j_2}$ are formed by consecutive points in the angular sorting, the overlap forces the starting indices to be close: $|j_2-j_1|\leq 7$. Therefore, for any fixed $i$ and $(p,q)\in D$, $\left|\{j:(i,j)\in f^{-1}(p,q)\}\right|\leq 8$. In other words, any $i$ can contribute at most 8 pairs $(i,j)$ to the preimage of $(p,q)$. Consequently, for any $(p,q)\in D$, $\left|\{i:(i,j)\in f^{-1}(p,q)\mbox{ for some $j$}\}\right|\geq \frac{1}{8}\left|f^{-1}(p,q)\right|$. But $v_i$ is visible from $(p,q)$ for all $i$ in $\{i:(i,j)\in f^{-1}(p,q)\mbox{ for some $j$}\}$. Therefore, $|V(p,q)|\geq \frac{1}{8}\left|f^{-1}(p,q)\right|$ for $(p,q)\in D$. Summing over all edges in $D$ yields:
\[\sum_{(p,q)\in D}|V(p,q)|\geq\tfrac{1}{8}\sum_{(p,q)\in D}\left|f^{-1}(p,q)\right|=\tfrac{1}{8}|I_2|.\]
		
Similarly, to obtain a bound relating the size of $\ps$ with the size of $|I_1|$, consider the function $g:I_1\rightarrow \ps$, defined by $g(i,j)=Q^i_j$. If $g(i_1,j_1)=g(i_2,j_2)=Q$, then the top point of $Q$ must be both $v_{i_1}$ and $v_{i_2}$, which implies $i_1=i_2$. Furthermore, since both $\spts^i_{j_1}$ and $\spts^i_{j_2}$ must contain the remaining four vertices $Q-\{v_i\}$, we have $|\spts^i_{j_1}\cap \spts^i_{j_2}|\geq 4$. As before, this implies $|j_2-j_1|\leq 5$ due to the angular sorting. Therefore,  $\left|g^{-1}(Q)\right|\leq 6$ for any $Q\in \ps$. Consequently, the size of $I_1$ is bounded by:
$|I_1|=\sum_{Q\in \ps}|g^{-1}(Q)|\leq 6|\ps|$, or
$|\ps|\geq \frac{1}{6}|I_1|$.

Adding the two inequalities, we obtain:
\[|\ps|+\sum_{(p,q)\in D}|V(p,q)|\geq \tfrac{1}{6}|I_1|+\tfrac{1}{8}|I_2|\geq \tfrac{1}{8}|I|=\tfrac{1}{8}\sum_{i=10}^n(i-9) = \tfrac{1}{16}(n^2-17n+72)\geq\tfrac{1}{16}(n^2-17n).\]
 Thus, $|\ps|+\sum_{(p,q)\in D}|V(p,q)|\geq\frac{1}{16}(n^2-17n)$.
		
By Corollary~\ref{vcp_corollary}, $5|\ps|+\sum_{(p,q)\in D}|C(p,q)|\geq \sum_{(p,q)\in D}|V(p,q)|$. Substituting this into our bound yields: \[6|\ps|+\sum_{(p,q)\in D}|C(p,q)|\geq |\ps|+\sum_{(p,q)\in D}|V(p,q)|\geq\tfrac{1}{16}(n^2-17n).\]
\end{proof}

\begin{figure}[tpb!]
	\centering\includegraphics[scale=.25]{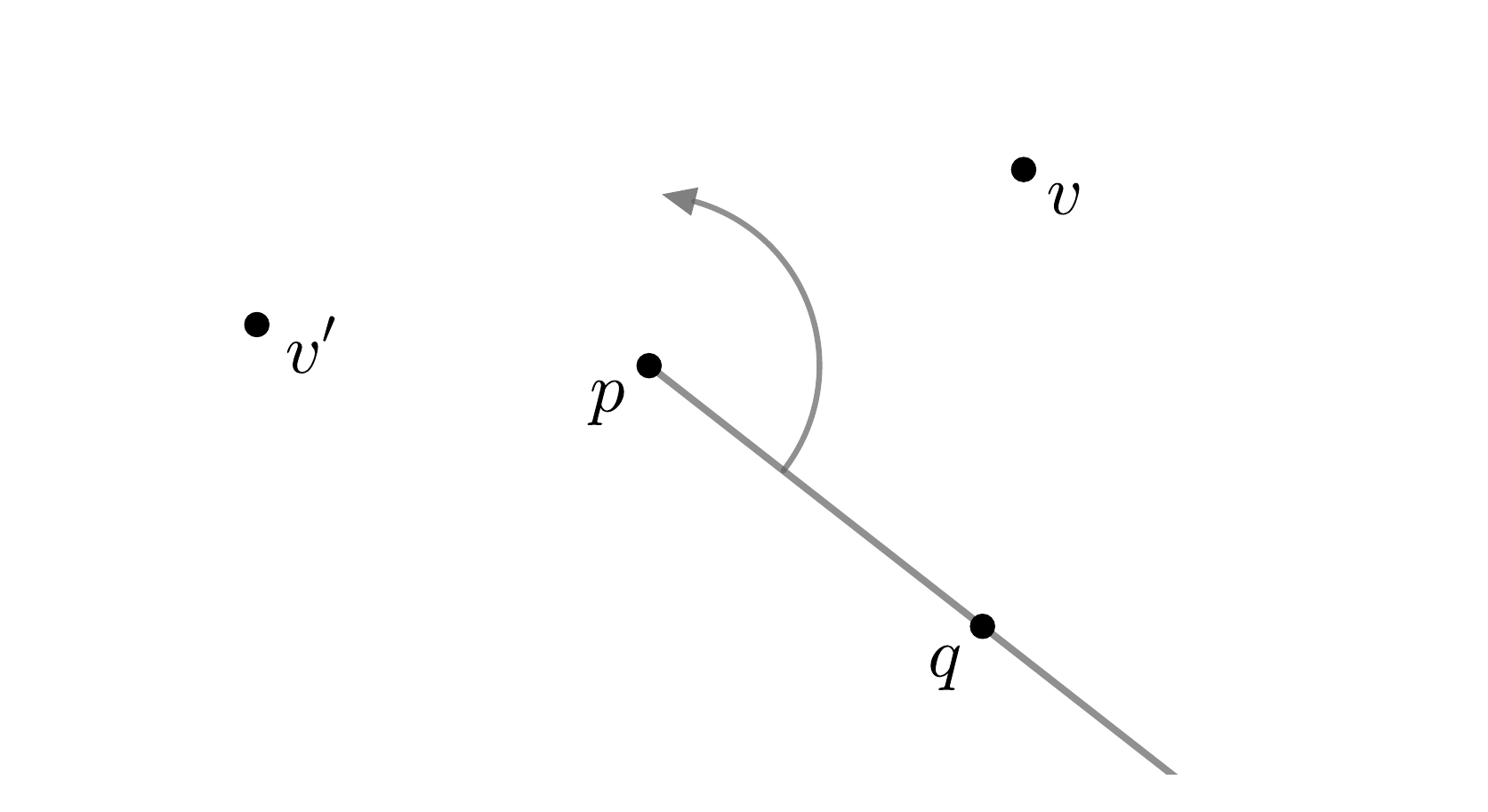}
	\caption{$v'\succ_{p,q}v$}
	\label{ccwOrder}
\end{figure}

\section{$C$-sets}
Now we consider the question of whether the lower bound on $\sum_{(p,q)\in D}|C(p,q)|$ from the previous theorem implies a lower bound on $\left|\bigcup_{(p,q)\in D}C(p,q)\right|$. In general, the answer is no, since $C$-sets may have large intersections (see, for example, Figure~\ref{commonCset}). However, for edges with a common left (or right) endpoint, $C$-sets form tree structures with relatively small intersections. Theorem~\ref{CCCT} formalizes this intersection structure, which we subsequently use in Theorem~\ref{boundCCC_theorem} to quantify the maximum overlap and establish a tight lower bound on the cardinality of the union. 

For any fixed $p$ and $q$ in $\spts$, define the counterclockwise order (ccw order for short) $\succ_{p,q}$ on $\spts-\{p,q\}$ with respect to $p$ and $q$ by $v'\succ_{p,q} v$ if  $q, v, v'$ appear in that order when sorted in counterclockwise order about $p$, for $v$ and $v'$ in $\spts-\{p,q\}$ and $v\neq v'$ (see Figure~\ref{ccwOrder}). Note that if $v\succeq v'$ and $v'\succeq v$ then $v=v'$ because no three points in $\spts$ are collinear. In general, we write $\succ$ instead of $\succ_{p,q}$ whenever $p$ and $q$ are clear from the context. Also, we extend $\succ_{p,q}$ to $\spts-\{p\}$ by allowing $v\succ_{p,q}q$ for all $v\in \spts-\{p,q\}$.

\begin{figure}[tpb!]
	\centering\includegraphics[scale=.6]{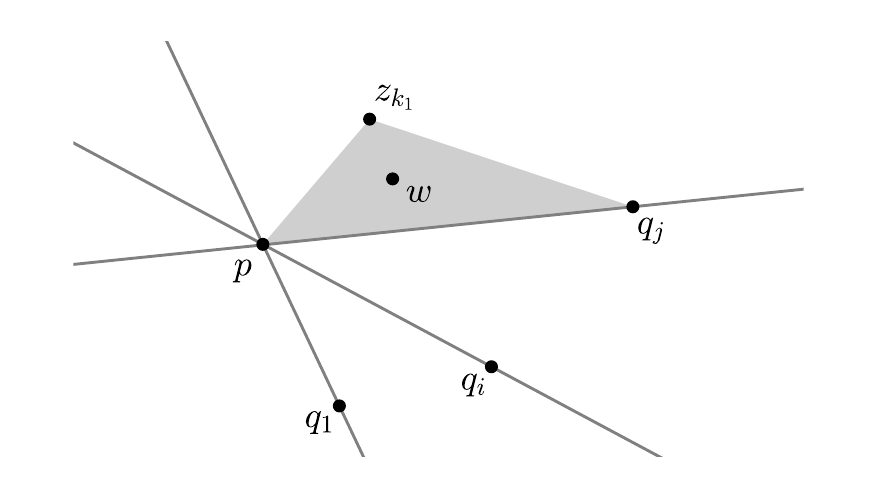}
	\caption{$z_{k_1}\succ w\succ q_j$}
	\label{CCCTfig_a}
\end{figure}

\begin{theorem} Let $h\geq 2$ and suppose $p, q_1, \ldots, q_h$, are vertices of $\spts$ such that $q_i\in\mc{H}^+(p,q_1)$ for all $2\leq i\leq h$ and $q_h\succ \cdots\succ q_1$, where $\succ$ is the ccw order with respect to $p$ and $q_1$. Let $1\leq j\leq h$ and suppose $C(p,q_j)=\{z_1,\ldots,z_{l}\}$, where $z_{l}\succ\cdots\succ z_1$ and $l\geq 0$. Then, 
	\begin{enumerate}
	\item For $1<j\leq h$, if $C(p,q_j)\cap\mc{H}^+(p,q_1)\neq\emptyset$ then there exists $k_1\geq 1$ such that: 
	\begin{itemize}
	\item $C(p,q_j)\cap\mc{H}^+(p,q_1)=\{z_1,\ldots,z_{k_1}\}$.
	\item For all $1\leq i<j$, if $w\in C(p,q_i)$ then either $w\succeq z_{k_1}$ or $q_j\succeq w$.
	\end{itemize}
	\item For $1\leq j< h$, if $C(p,q_j)\cap\mc{H}^+(p,q_h)\neq\emptyset$ then there exists $k_2\geq 1$ such that: 
	\begin{itemize}
		\item $C(p,q_j)\cap\mc{H}^+(p,q_h)=\{z_{k_2},\ldots,z_l\}$.
		\item For all $j<i\leq h$, if $w\in C(p,q_i)$ then either $z_{k_2}\succeq w$ or $w\in\mc{H}^-(p,q_j)$.
	\end{itemize}
	\end{enumerate}
\label{CCCT}
\end{theorem}

\begin{figure}[tpb!]
	\centering\includegraphics[scale=.55]{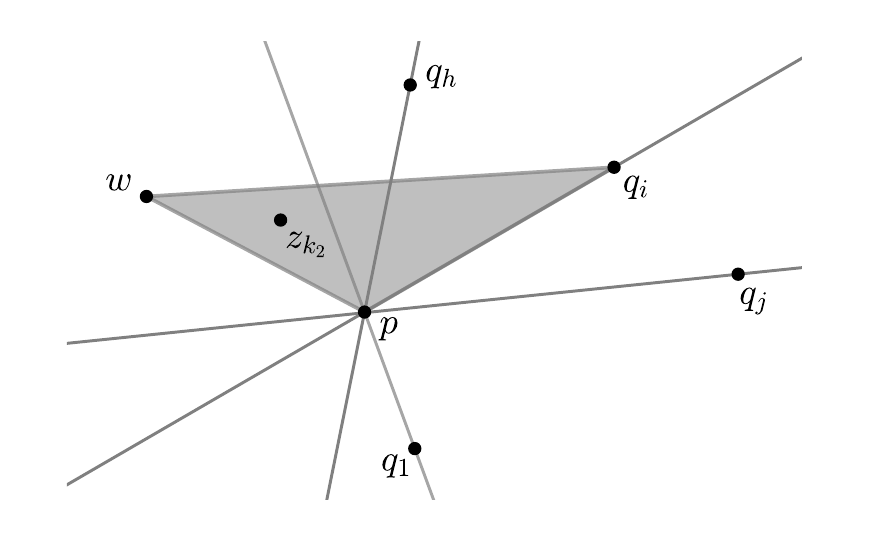}
	\caption{$w\succ z_{k_2}\succ q_i$}
	\label{CCCTfig_b}
\end{figure}

\begin{proof} Let $\succ$ be the ccw order with respect to $p$ and $q_1$. Suppose $h\geq 2$, $q_i\in\mc{H}^+(p,q_1)$ for all $2\leq i\leq h$, and $q_h\succ \cdots\succ q_2$. For the sake of brevity, we write $C_i$ for $C(p,q_i)$, $C^+_i$ for $C(p,q_i)\cap\mc{H}^+(p,q_1)$ and $C^-_i$ for $C(p,q_i)\cap\mc{H}^+(p,q_h)$.

Let $1\leq j\leq h$ and suppose $C_j=\{z_1,\ldots,z_{l}\}$, where $z_{l}\succ\cdots\succ z_1$. Note that the  order $z_{l}\succ\cdots\succ z_1$ is the order in which $z_1,\ldots,z_{l}$ appear on the boundary of  $\conv(\spts\cap\mc{H}^+(p,q_j))$ if we traverse this boundary in the clockwise direction starting at $z_1$. 

To prove the first part of the theorem, 
suppose $1<j\leq h$ and $C_j^+\neq\emptyset$. Let $k_1$ be the largest $k$ such that $z_k\in\mc{H}^+(p,q_1)$. Clearly, $z_k\in C_j^+$ for all $k\leq k_1$ and $z_k\not\in C_j^+$ for all $k>k_1$ since $z_{l}\succ\cdots\succ z_1$. Therefore, $C_j^+=\{z_1,\ldots,z_{k_1}\}$. 

Now suppose $1\leq i<j$ and $w\in C_i$. In order to obtain a contradiction, suppose $z_{k_1}\succ w\succ q_j$ (see Figure~\ref{CCCTfig_a}). Note that $q_j$ belongs to $\mc{H}^+(p,q_i)$ since both $q_i$ and $q_j$ belong to $\mc{H}^+(p,q_1)$ and $q_j\succ q_i$. For the same reason, $C_j^+\subset\mc{H}^+(p,q_j)\cap\mc{H}^+(p,q_1)\subset \mc{H}^+(p,q_i)\cap\mc{H}^+(p,q_1)$, so $z_{k_1}$ belongs to $\mc{H}^+(p,q_i)$ too. Since $w\in C_i$, there exists a line through $w$ that separates $\{z_{k_1},q_j\}$ (and the rest of $\spts\cap\mc{H}^+(p,q_i)-\{w\}$) from $\{p,q_i\}$. In particular, this line separates $\{z_{k_1},q_j\}$ from $p$. Therefore, $w$ belongs to the interior of $\conv(p,q_j,z_{k_1})$. This is a contradiction because $\conv(p,q_j,z_{k_1})$ is contained in the region between $\mc{L}(p,q_j)$ and $C_j$ (as $z_{k_1}\in C_j$), which is empty. Hence, $z_{k_1}\succ w\succ q_j$ must be false. Therefore, either $w\succeq z_{k_1}$ or $q_j\succeq w$.

To prove the second part of the theorem, suppose $1\leq j< h$ and $C_j^-\neq\emptyset$. Let $k_2$ be the smallest $k$ such that $z_k\in\mc{H}^+(p,q_h)$. Clearly, $z_k\in C_j^-$ for all $k\geq k_2$ and $z_k\not\in C_j^-$ for all $k<k_2$ since $z_{l}\succ\cdots\succ z_1$. Hence, $C_j^-=\{z_{k_2},\ldots,z_l\}$.

Suppose $j<i\leq h$ and $w\in C_i$. Additionally, suppose $w\in\mc{H}^+(p,q_j)$. In order to obtain a contradiction, suppose $w\succ z_{k_2}$ (see Figure~\ref{CCCTfig_b}). Clearly, $w\succ z_{k_2}\succ q_i$ since $z_{k_2}\succ q_h\succeq q_i$. Note that $q_i$ belongs to $\mc{H}^+(p,q_j)$ since both $q_i$ and $q_j$ belong to $\mc{H}^+(p,q_1)$ and $q_i\succ q_j$. Hence, both $w$ and $q_i$ belong to $\mc{H}^+(p,q_j)$. Since $z_{k_2}\in C_j$, there exists a line through $z_{k_2}$ that separates $\{w,q_i\}$ (and the	 rest of $\spts\cap\mc{H}^+(p,q_j)-\{z_{k_2}\}$) from $\{p,q_j\}$. In particular, this line separates $\{w,q_i\}$ from $p$. Therefore, $z_{k_2}$ necessarily belongs to the interior of $\conv(p,q_i,w)$. This is a contradiction because $\conv(p,q_i,w)$ is contained in the region between $\mc{L}(p,q_i)$ and $C_i$ (as $w\in C_i$), which is empty. Hence, $w\succ z_{k_2}$ must be false. Therefore, $z_{k_2}\succeq w$ when $w\in\mc{H}^+(p,q_j)$. If $w\not\in\mc{H}^+(p,q_j)$, note that $w\neq q_j$ as $w\in\mc{H}^+(p,q_i)$ but $q_j\in\mc{H}^-(p,q_i)$ since $i>j$. Hence, $w\in\mc{H}^-(p,q_j)$ when $w\not\in\mc{H}^+(p,q_j)$. Therefore, for all $j<i\leq h$ and $w\in C_i$, either  $z_{k_2}\succeq w$ or $w\in\mc{H}^-(p,q_j)$.  
\end{proof}

\begin{theorem} Let $h\geq 2$ and suppose $p, q_1, \ldots, q_h$ are distinct vertices of $\spts$ such that $q_i\in\mc{H}^+(p,q_1)$ for all $2\leq i\leq h$. Then, $\left|\bigcup_{j=1}^h C(p,q_j)\right|\geq \left(\sum_{j=1}^h\left|C(p,q_j)\right|\right)-2h+3$.
	\label{boundCCC_theorem}
\end{theorem}

\begin{figure}[tpb!]
	\centering\includegraphics[scale=.325]{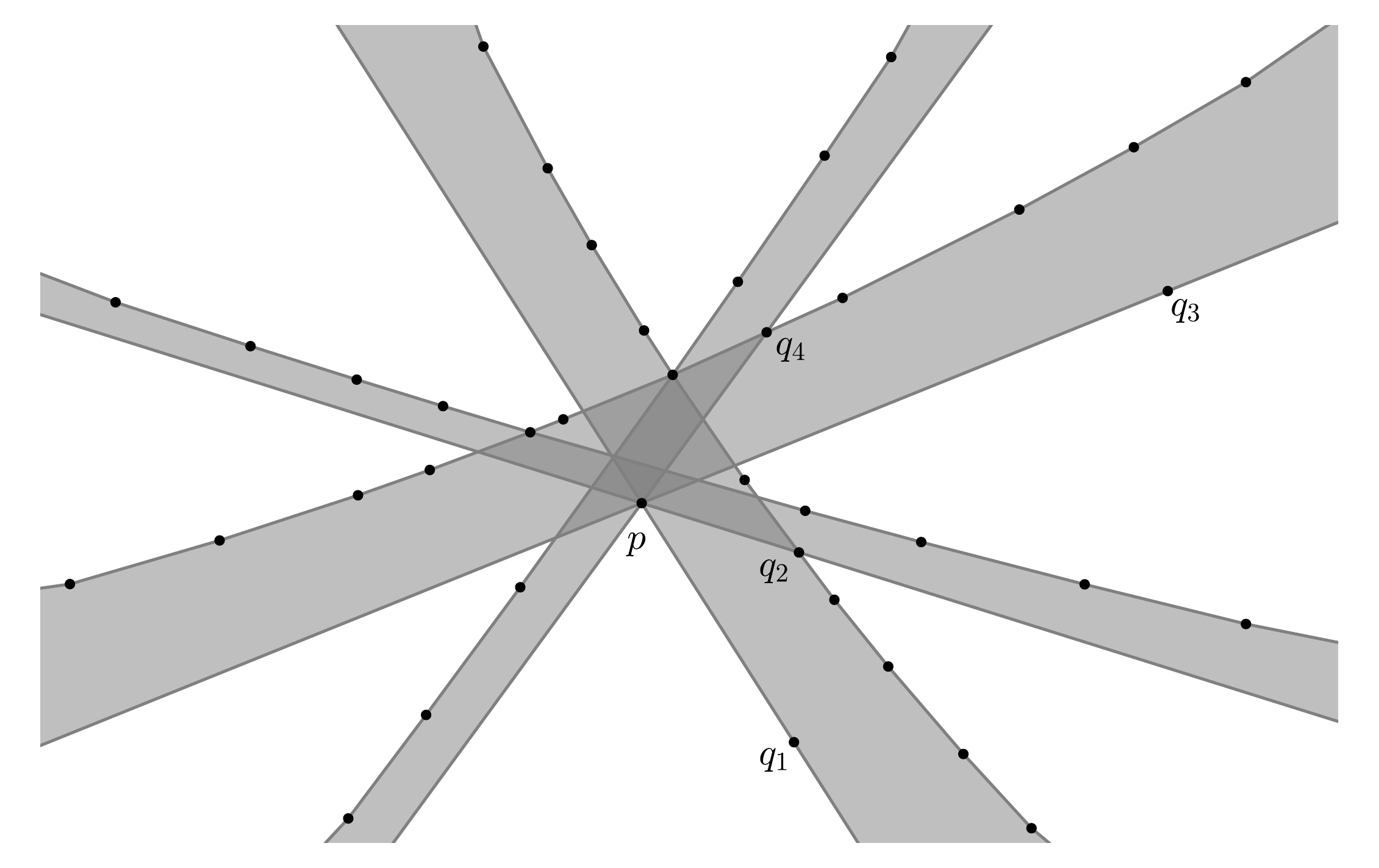}
	\caption{$\left|C_3^+\cap(C_1\cup C_2)\right|=1$ and $\left|C_3^-\cap C_4\right|=1$.}
	\label{Figure9b}
\end{figure}

\begin{proof} Suppose $h\geq 2$ and let $p, q_1, \ldots, q_h$ be distinct vertices of $\spts$ such that $q_i\in\mc{H}^+(p,q_1)$ for all $2\leq i\leq h$. Let $\succ$ be the ccw order with respect to $p$ and $q_1$. Without loss of generality, assume $q_h\succ\cdots\succ q_1$. Let $C_i$ denote $C(p,q_i)$, $C_i^+$ denote $C(p,q_i)\cap\mc{H}^+(p,q_1)$, and $C^-_i$ denote $C(p,q_i)\cap\mc{H}^+(p,q_h)$. See Figure~\ref{Figure9b} for an example. 
		
First, we consider the restriction of the $C$-sets to $\mc{H}^+(p,q_1)$. 
Let $1<j\leq h$ and suppose $C_j^+\neq\emptyset$. By Theorem~\ref{CCCT}, there exists $z'\in C_j^+$ such that $w\succeq z'$ or $q_j\succeq w$ for all $1\leq i<j$ and $w\in C_i$. Moreover, $z'\succeq z$ for all $z\in C_j^+$. Therefore, $w\succeq z'\succeq z$ or $q_j\succeq w$ for all $w\in C_i$, $1\leq i<j$, and $z\in C_j^+$. Hence, if $w=z$ for some $w\in C_i$ (with $1\leq i<j$) and $z\in C_j^+$, then necessarily $w=z'=z$. Therefore, $C_j^+\cap\left(\bigcup_{i=1}^{j-1}C_i\right)\subseteq\{z'\}$. In particular, $C_j^+\cap\left(\bigcup_{i=1}^{j-1}C_i^+\right)\subseteq\{z'\}$. 
	
Hence, $\left| C_j^+\cap\left(\bigcup_{i=1}^{j-1}C_i^+\right)\right|\leq 1$, for $1<j\leq h$. Since the sets $C_j^+-\bigcup_{i=1}^{j-1}C_i^+$ are disjoint and their union is $\bigcup_{j=1}^hC^+_j$, for $1\leq j\leq h$, $\left|\bigcup_{j=1}^hC^+_j\right|=\left|\bigcup_{j=1}^h\left(C^+_j-\bigcup_{i=1}^{j-1}C^+_i\right)\right|=\sum_{j=1}^h\left|C^+_j-\bigcup_{i=1}^{j-1}C^+_i\right|$.
On the other hand, $\left|C^+_j-\bigcup_{i=1}^{j-1}C^+_i\right|=|C^+_j|-\left|C^+_j\cap\bigcup_{i=1}^{j-1}C^+_i\right|\geq|C^+_j|-1$, for $1<j\leq h$, which implies $\sum_{j=1}^h\left|C^+_j-\bigcup_{i=1}^{j-1}C^+_i\right|\geq|C^+_1|+\sum_{j=2}^h(|C^+_j|-1)=\left(\sum_{j=1}^h|C^+_j|\right)-h+1$. Therefore, $\left|\bigcup_{j=1}^hC^+_j\right|\geq\left(\sum_{j=1}^h|C^+_j|\right)-h+1$. 
	
Now we consider the restriction of the $C$-sets to $\mc{H}^+(p,q_h)$. Let $1\leq j<h$ and suppose $C_j^-\neq\emptyset$. By Theorem~\ref{CCCT}, there exists $z''$ in $C_j^-$ such that $z''\succeq w$ or $w\in\mc{H}^-(p,q_j)$ for all $j<i\leq h$ and $w\in C_i$. Additionally, $z\succeq z''$ for all $z\in C_j^-$. Hence, $z\succeq z''\succeq w$ or $w\in\mc{H}^-(p,q_j)$ for all $w\in C_i$, $j<i\leq h$, and $z\in C_j^-$. Therefore, if $z=w$ for $w\in C_i$ (with $j<i\leq h$) and $z\in C_j^-$, then necessarily $z=z''=w$. Hence, $C_j^-\cap\left(\bigcup_{i=j+1}^{h}C_i\right)\subseteq\{z''\}$. In particular, $C_j^-\cap\left(\bigcup_{i=j+1}^{h}C_i^-\right)\subseteq\{z''\}$. 
Note that $\mc{H}^+(p,q_j)\subset\mc{H}^+(p,q_1)\cup\mc{H}^+(p,q_h)$, which implies 
$C_j\cap\mc{H}^-(p,q_1)\subseteq C_j\cap\mc{H}^+(p,q_h)=C_j^-$. To avoid overcounting, we need to consider only the part of $C_j^-$ that is contained in $\mc{H}^-(p,q_1)$. Let $C_j^=$ denote $C_j\cap\mc{H}^-(p,q_1)$. Clearly, $C_j^{=}=C_j^-\cap\mc{H}^-(p,q_1)$. Then, $C_j^=\cap\left(\bigcup_{i=j+1}^{h}C_i^=\right)=\left(C_j^- \cap\mc{H}^-(p,q_1)\right)\cap\left(\bigcup_{i=j+1}^{h}\left(C_i^-\cap\mc{H}^-(p,q_1)\right)\right)=C_j^-\cap\left(\bigcup_{i=j+1}^{h}C_i^-\right)\cap\mc{H}^-(p,q_1)\subseteq\{z''\}\cap\mc{H}^-(p,q_1)\subseteq\{z''\}$.
	
Therefore, $\left|C^=_j\cap\left(\bigcup_{i=j+1}^hC_i^=\right)\right|\leq 1$, for $1\leq j<h$. Since $C_1^{=}=\emptyset$, we only need to consider $j \in \{2, \ldots, h\}$. Note that, for $j\in\{2,\ldots,h\}$, the sets $C_j^=-\bigcup_{i=j+1}^{h}C_i^=$ are disjoint and their union is $\bigcup_{j=2}^hC^=_j$. Hence, $\left|\bigcup_{j=2}^hC^=_j\right|=\left|\bigcup_{j=2}^h\left(C^=_j-\bigcup_{i=j+1}^hC^=_i\right)\right|=\sum_{j=2}^h\left|C^=_j-\bigcup_{i=j+1}^hC^=_i\right|$.
On the other hand, $\left|C^=_j-\bigcup_{i=j+1}^hC^=_i\right|=|C^=_j|-\left|C^=_j\cap\bigcup_{i=j+1}^hC^=_i\right|\geq|C^=_j|-1$, for $2\leq j<h$, which implies $\sum_{j=2}^h\left|C^=_j-\bigcup_{i=j+1}^hC^=_i\right|\geq|C^=_h|+\sum_{j=2}^{h-1}(|C^=_j|-1)=\left(\sum_{j=2}^h|C^=_j|\right)-h+2$. Therefore, $\left|\bigcup_{j=2}^hC^=_j\right|\geq\left(\sum_{j=2}^h|C^=_j|\right)-h+2$.

Since $C_j=C_j^{=}\cup C_j^+$ for all $1\leq j\leq h$ and $C_{j_1}^{=}\cap C_{j_2}^+=\emptyset$ for all $j_1$ and $j_2$, we have:
$\left|\bigcup_{j=1}^h C_j\right| = \left|\bigcup_{j=1}^h \left(C_j^{=}\cup C_j^+\right)\right| = \left|\bigcup_{j=2}^h C_j^{=}\right|+\left|\bigcup_{j=1}^h C^+_j\right| \geq \left(\sum_{j=2}^h\left|C^=_j\right|\right)-h+2+ \left(\sum_{j=1}^h\left|C^+_j\right|\right)-h+1 \geq \left(\sum_{j=1}^h\left|C_j\right|\right)-2h+3$, 
which completes the proof. 
\end{proof}

\begin{figure}[tpb!]
	\centering\includegraphics[scale=.375]{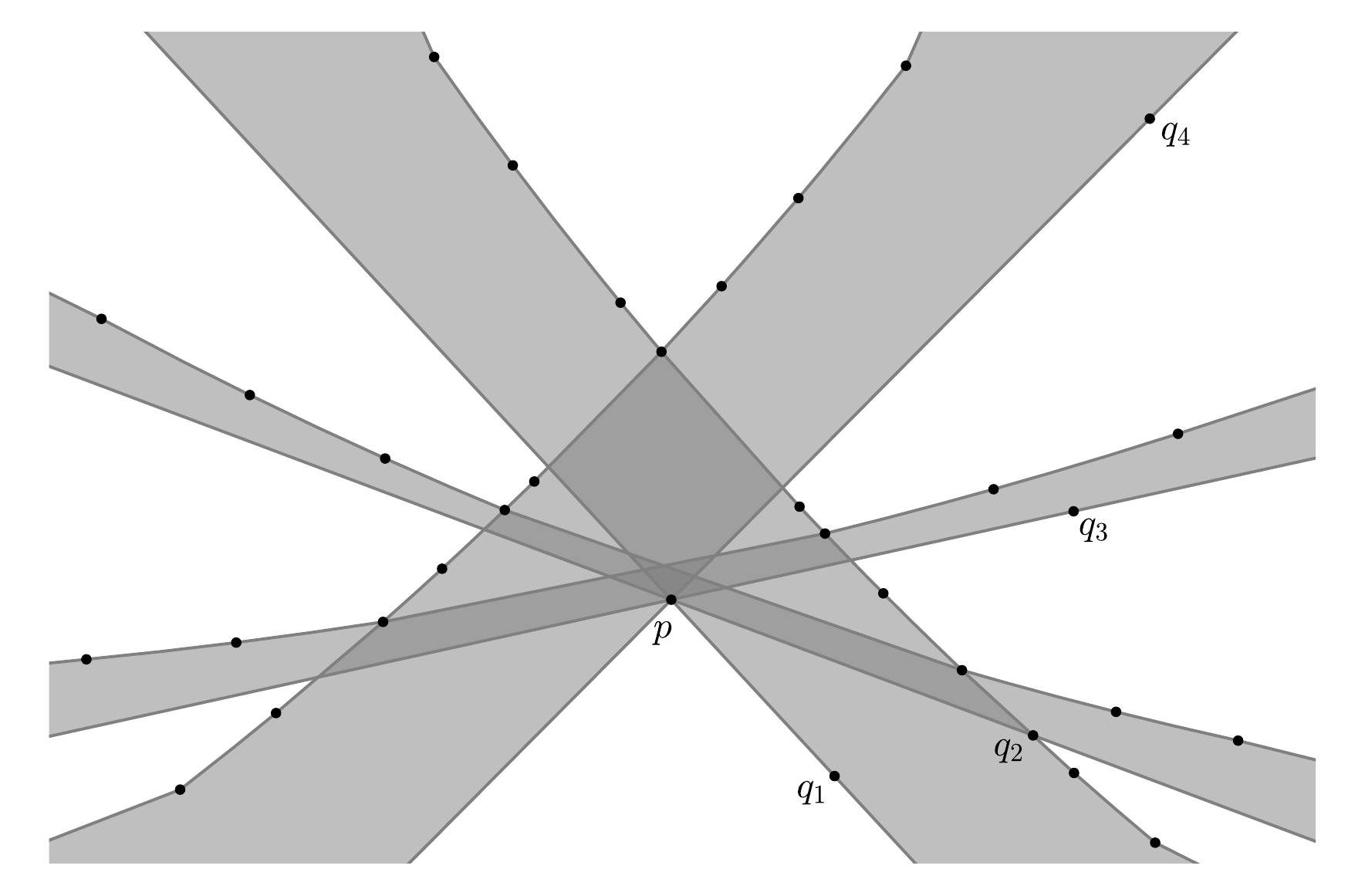}
	\caption{$\left|\bigcup_{j=1}^4 C_j\right|= \left(\sum_{j=1}^4\left|C_j\right|\right)-2(4)+3$.}
	\label{Figure10b}
\end{figure}

Note that the bound in Theorem~\ref{boundCCC_theorem} is tight (see Figure~\ref{Figure10b} for an example with $h=4$).

\begin{corollary} Let $p, q_1, \ldots, q_h$, $h\geq 1$, be distinct vertices of $\spts$. Then, $\sum_{i=1}^h |C(p,q_i)|\leq 4(n-1)$.
\label{upperboundCCC_corollary}
\end{corollary}

\begin{proof} We repeatedly use the following obvious facts: $\left|\bigcup_{j} C(p,q_j)\right|\leq |\spts-\{p\}|= n-1$ and $h\leq n-1$. If $h=1$, the conclusion is immediate since $|C(p,q_1)|\leq n-1 \leq 4(n-1)$. 
If $h\geq 2$, assume without loss of generality that $q_h\succ\cdots\succ q_1$, where $\succ$ is the ccw order with respect to $p$ and $q_1$. Let $C_j$ denote $C(p,q_j)$. 
	
If $q_j\in\mc{H}^+(p,q_1)$ for all $2\leq j\leq h$, then by Theorem~\ref{boundCCC_theorem} we have:
$\sum_{j=1}^h|C_j|\leq\left|\bigcup_{j=1}^h C_j\right|+2h-3\leq (n-1)+2h-3\leq 4(n-1)$. 
If not, let $k\geq 1$ be the minimum $j$ such that $q_{j+1}\in\mc{H}^-(p,q_1)$. Then, $q_j\in\mc{H}^+(p,q_1)$ for $2\leq j\leq k$, and $q_j\in\mc{H}^-(p,q_1)$ for $k+1\leq j\leq h$. 
	
We bound the sum in two parts: $\sum_{j=1}^k|C_j|$ and $\sum_{j=k+1}^h|C_j|$. 
For the first part, if $k=1$, we trivially have $\sum_{j=1}^1|C_j| = |C_1| \leq n-1 = n+2k-3$. If $k\geq 2$, applying Theorem~\ref{boundCCC_theorem} to $p, q_1,\ldots, q_k$ yields:
$\sum_{j=1}^k|C_j|\leq\left|\bigcup_{j=1}^k C_j\right|+2k-3\leq (n-1)+2k-3=n+2k-4\leq n+2k-3.$
	
For the second part, if $k+1=h$, then $h-k=1$, and we trivially have $\sum_{j=k+1}^h|C_j| = |C_h| \leq n-1 = n+2(h-k)-3$. If $k+1<h$, applying Theorem~\ref{boundCCC_theorem} to $p, q_{k+1},\ldots, q_h$ (since $q_j\in\mc{H}^+(p,q_{k+1})$ for all $k+1< j\leq h$) yields:
$\sum_{j=k+1}^h|C_j|\leq\left|\bigcup_{j=k+1}^h C_j\right|+2(h-k)-3\leq (n-1)+2(h-k)-3=n+2(h-k)-4 \leq n+2(h-k)-3.$
	
Therefore, combining the two parts, we obtain:
$\sum_{j=1}^h|C_j|=\sum_{j=1}^k|C_j|+\sum_{j=k+1}^h|C_j|\leq (n+2k-3)+(n+2(h-k)-3)=2n+2h-6\leq 2n+2(n-1)-6 = 4n-8 \leq 4(n-1)$.
\end{proof}

By grouping all sets $C(p,q)$ by their left endpoint $p$, applying Corollary~\ref{upperboundCCC_corollary}, and dividing by the total number of oriented edges, $n(n-1)$, we find that the average size of a $C$-set in any $n$-point set $\spts$ is bounded above by $4$. This is undoubtedly a loose bound. We think that improving it is an interesting direction for future research.

For any points $p$ and $q$ in $\spts$ with $p\neq q$, and points $a$, $b$, and $q'$ in $\spts$ such that $b\succeq_{p,q} a\succ_{p,q} q'$, we define the subinterval $[b,a]_{q'}$ of $C(p,q')$ to be $\{z\in C(p,q'): b\succeq z\succeq a\}$. Clearly, $[b,a]_{q'}=C(p,q')\cap\left(\{a,b\}\cup\left(\mc{H}^+(p,a)\cap \mc{H}^-(p,b)\right)\right)$.  
Similarly, we define $(b,a]_{q'}$ to be $\{z\in C(p,q'): b\succ z\succeq a\}$. Note that $(b,a]_{q'}$ equals $[b,a]_{q'}$ if and only if $b$ is not in $C(p,q')$. Note that intervals may be empty, and the points $p$ and $q$ are left to be inferred from the context in this notation. 

Given two subsets $U$ and $U'$ of $\spts-\{p,q\}$, we write $U'\succ_{p,q} U$ if $v'\succ_{p,q} v$ for all $v'\in U'$ and $v\in U$. Also, we say that $U$ and $U'$ are angle-disjoint under $\succ_{p,q}$ if $U'\succ U$ or $U\succ U'$. A family of sets is angle-disjoint under $\succ_{p,q}$ if any two of its elements are angle-disjoint, in which case $\succ_{p,q}$ induces a total order on the family. Note that the empty set is angle-disjoint from any set. 

The next theorem shows that, under certain conditions, it is possible to partition $\bigcup_{i=1}^h C(p,q_i)$ into a family of at most $4h-2$ angle-disjoint subintervals of $C(p,q_1),\ldots,C(p,q_h)$. 

\begin{theorem} Let $h\geq 2$ and suppose $p,q_1,\ldots,q_h$ are distinct vertices of $\spts$ such that $q_i\in\mc{H}^+(p,q_1)$ for all $2\leq i\leq h$. Let $\succ$ be the ccw order with respect to $p$ and $q_1$. Then, there exists a set $J$ such that 
\begin{itemize}
	\item For all $X\in J$, $X\neq\emptyset$ and there exist $1\leq i\leq h$ and $\{a,b\}\subseteq C(p,q_i)$ such that $X=[b,a]_{q_i}$ or $X=(b,a]_{q_i}$. In other words, $J$ is a family of nonempty subintervals of $C(p,q_1), \ldots, C(p,q_h)$. 
	\item $\bigcup_{X\in J} X=\bigcup_{i=1}^h C(p,q_i)$.
	\item If $\{X,X'\}\subseteq J$ and $X\neq X'$ then $X\succ X'$ or $X'\succ X$, i.e., $X$ and $X'$ are angle-disjoint. 
	\item $|J|\leq 4h-2$.
\end{itemize}
\label{CCCpartition_theorem}
\end{theorem}

\begin{proof} Let $h\geq 2$ and suppose $p,q_1,\ldots,q_h$ are distinct vertices of $\spts$ such that $q_i\in\mc{H}^+(p,q_1)$ for all $2\leq i\leq h$. Let $\succ$ be the ccw order with respect to $p$ and $q_1$. Without loss of generality, assume $q_h\succ\cdots\succ q_1$. Let $C_i$ denote $C(p,q_i)$, $C_i^+$ denote $C(p,q_i)\cap\mc{H}^+(p,q_1)$ and $C^-_i$ denote $C(p,q_i)\cap\mc{H}^+(p,q_h)$.
	
First, we show that there exists a partition $J_h^+$ of $\bigcup_{i=1}^{h} C_i^+$ into no more than $2h-1$ angle-disjoint subintervals of $C_1^+,\ldots,C_{h}^+$. Proceeding inductively, we will show that for every $1\leq h'\leq h$ there exists $J_{h'}^+$ satisfying the following conditions: (1) every element of $J_{h'}^+$ is a subinterval of $C_i^+$ for some $1\leq i\leq h'$, (2) $\bigcup_{X\in J_{h'}^+}X=\bigcup_{i=1}^{h'}C_i^+$, (3) $X\succ X'$ or $X'\succ X$ for all $\{X,X'\}\subseteq J_{h'}^+$ with $X\neq X'$, and (4) $\left|J_{h'}^+\right|\leq 2h'-1$. Let $J_1^+=\{C_1^+\}$. Note that $C_1^+\neq\emptyset$ because $q_2\in\mc{H}^+(p,q_1)$. Clearly, $J_1^+$ satisfies the four conditions for $h'=1$. Now suppose $J_{j}^+$ satisfies the four conditions above for $h'=j$, where $1\leq j<h$. If $C_{j+1}^+=\emptyset$ then $J_{j+1}^+=J_{j}^+$ satisfies the conditions for $h'=j+1$. If $C_{j+1}^+\neq\emptyset$, let $C_{j+1}=\{z_1,\ldots,z_l\}$, with $z_l\succ\cdots\succ z_1$. By the first part of Theorem~\ref{CCCT}, there exists $k_1\geq 1$ such that $C_{j+1}^+=\{z_1,\ldots,z_{k_1}\}=[z_{k_1},z_1]_{q_{j+1}}$ and, for all $1\leq i<j+1$, if $w\in C_i$ then either $w\succeq z_{k_1}$ or $q_{j+1}\succeq w$. 
	
Suppose $X\succeq \{z_{k_1}\}$ or $\{q_{j+1}\}\succeq X$ for all $X\in J_j^+$. Then, $X\succeq [z_{k_1},z_1]_{q_{j+1}}$ or $[z_{k_1},z_1]_{q_{j+1}}\succ X$, for all $X\in J_j^+$. We consider two cases. If $z_{k_1}\in X'$ for some $X'\in J_j^+$, then $z_{k_1}$ is already covered by $J_j^+$. Since $X\succeq [z_{k_1},z_1]_{q_{j+1}}$ or $[z_{k_1},z_1]_{q_{j+1}}\succ X$ implies $X\succ (z_{k_1},z_1]_{q_{j+1}}$ or $(z_{k_1},z_1]_{q_{j+1}}\succ X$, it follows that $J_{j+1}^+=J_j^+\cup\{(z_{k_1},z_1]_{q_{j+1}}\}$ is a partition of $\bigcup_{i=1}^{h'}C_i^+$ satisfying the properties above for $h'=j+1$. If not, i.e., if $z_{k_1}\not\in X'$ for all $X'\in J_j^+$, then we have $X\succ [z_{k_1},z_1]_{q_{j+1}}$ or $[z_{k_1},z_1]_{q_{j+1}}\succ X$ for all $X\in J_j^+$, so $J_{j+1}^+=J_j^+\cup\{[z_{k_1},z_1]_{q_{j+1}}\}$ is the desired partition of $\bigcup_{i=1}^{h'}C_i^+$ for $h'=j+1$. 
	
Now suppose there exists $X_0\in J_j^+$ such that $X_0\not\succeq\{z_{k_1}\}$ and $\{q_{j+1}\}\not\succeq X_0$. Then, there exist $w_1$ and $w_2$ in $X_0$ such that $z_{k_1}\succ w_1$ and $w_2\succ q_{j+1}$. This implies $q_{j+1}\succeq w_1$ and $w_2\succeq z_{k_1}$, by the Theorem~\ref{CCCT} inequalities above, since $\{w_1,w_2\}\subseteq X_0\subseteq C_i$ for some $1\leq i<j+1$. Hence, if $X_0=[b,a]_{q_i}$ then necessarily $b\succeq z_{k_1}$ and $q_{j+1}\succeq a$. Let $X_0^2=\{w\in X_0 : w\succeq z_{k_1}\}$ and $X_0^1=\{w\in X_0 : q_{j+1}\succeq w\}$. Clearly, $X_0^1$ and $X_0^2$ are nonempty intervals satisfying $X_0=X_0^1\cup X_0^2$ and $X_0^2\succeq [z_{k_1},z_1]_{q_{j+1}}\succ X_0^1$. Furthermore, for all $X\in J_{j}^+$ with $X\neq X_0$, either $X\succ X_0$ or $X_0\succ X$ by the properties of $J_j^+$. Therefore, either $X\succ X_0^2\succeq [z_{k_1},z_1]_{q_{j+1}}\succ X_0^1$ or $X_0^2\succeq [z_{k_1},z_1]_{q_{j+1}}\succ X_0^1\succ X$, for all $X\in J_{j}^+$ with $X\neq X_0$. As before, there are two cases. If $z_{k_1}\in X_0^2$ then $z_{k_1}$ is already covered in $J_j^+$ and we have $X_0^2\succ (z_{k_1},z_1]_{q_{j+1}}\succ X_0^1$. Let $J_{j+1}^+=(J_j^+-\{X_0\})\cup(\{X_0^1,X_0^2,(z_{k_1},z_1]_{q_{j+1}}\}-\{\emptyset\})$. Clearly, $|J_{j+1}^+|\leq|J_j^+|-1+3\leq (2j-1)-1+3=2(j+1)-1$. Therefore, $J_{j+1}^+$ satisfies the required properties for $h'=j+1$. For the second case, assume $z_{k_1}\not\in X_0^2$. Then, $X_0^2\succ [z_{k_1},z_1]_{q_{j+1}}\succ X_0^1$ so $J_{j+1}^+=(J_j^+-\{X_0\})\cup\{X_0^1,X_0^2,[z_{k_1},z_1]_{q_{j+1}}\}$ is the desired partition. By induction, it follows that a partition exists for all $1\leq h'\leq h$.
	
A similar inductive argument, but based on the second part of Theorem~\ref{CCCT} and using $C_h^-$ as the base case, shows that there exists a partition $J_h^-$ of $\bigcup_{i=1}^{h} C_i^-$ into at most $2h-1$ angle-disjoint subintervals of $C_1^-,\ldots,C_{h}^-$.
	
For any set $X$, let $X^{=}$ denote $X\cap\mc{H}^-(p,q_1)$. Now let $J_h^{=}=\{X^{=}: X\in J_h^-\}-\{\emptyset\}$. Clearly, $J_h^{=}$ is a partition of $\bigcup_{i=1}^h C_i^{=}$ into at most $2h-1$ angle-disjoint subintervals of $C_1^{=},\ldots,C_h^{=}$ because the properties of $J_h^-$ are inherited by $J_h^{=}$ (in particular, the intersection of a subinterval of $C_i$ with $\mc{H}^-(p,q_1)$ is a subinterval of $C_i^{=}$).  
	
Finally, note that $X'\succ X$ for all $X'\in J_h^{=}$ and $X\in J_h^+$ because $X'\subset\mc{H}^-(p,q_1)$ and $X\subset\mc{H}^+(p,q_1)$. Therefore, $J_h=J_h^{=}\cup J_h^+$ is a partition of $\bigcup_{i=1}^h C_i$ into at most $4h-2$ angle-disjoint subintervals of $C_1,\ldots,C_h$, since
$\left(\bigcup_{i=1}^h C_i^{=}\right)\cup\left(\bigcup_{i=1}^h C_i^{+}\right)=\bigcup_{i=1}^h \left(C_i^{=}\cup C_i^+\right)=\bigcup_{i=1}^h C_i$.
\end{proof}

In the bound $4h-2$ of our previous theorem, the coefficient of $h$ is the best possible, but the constant term is not.

We remark that the linear complexity of the union of $C$-sets proven above is an instance of a well-known result concerning the complexity of the union of regions whose boundaries intersect a constant number of times and share a common point (see~\cite{AgarwalSharir1995}), as is the case for the empty regions between the edges and their corresponding extended polygonal chains.

While it is interesting to know that the sets in the previous theorem are angle-disjoint subintervals of $C$-sets, for the purpose of counting pentagons, we only need the fact that they are in convex position and that their convex hulls are pairwise disjoint.

\begin{corollary} Let $h\geq 2$ and suppose $p,q_1,\ldots,q_h$ are distinct vertices of $\spts$ such that $q_i\in\mc{H}^+(p,q_1)$ for all $2\leq i\leq h$. Then, there exists a set $J$ such that:
	\begin{itemize}
		\item For all $X\in J$, $X$ is a nonempty subset of $\spts$ in convex position.
		\item If $\{X, X'\}\subseteq J$ and $X\neq X'$, then $\mbox{\emph{conv}}(X)\cap\mbox{\emph{conv}}(X')=\emptyset$.
		\item $|J|\leq 4h-2$.
		\item $\sum_{X\in J} |X|\geq \left(\sum_{i=1}^h |C(p,q_i)|\right)-2h+3$.
	\end{itemize} 
	\label{convex1_corollary}
\end{corollary}

\begin{proof} By Theorem~\ref{CCCpartition_theorem}, there exists a partition $J$ of $\bigcup_{i=1}^h C(p,q_i)$ into at most $4h-2$ angle-disjoint subintervals of $C(p,q_1),\ldots,C(p,q_h)$. The first three properties in the conclusion follow immediately. For the fourth, because $J$ is a partition of $\bigcup_{i=1}^h C(p,q_i)$, we have:
$\sum_{X\in J}|X|=\left|\bigcup_{X\in J}X\right|=\left|\bigcup_{i=1}^h C(p,q_i)\right|$. 
By Theorem~\ref{boundCCC_theorem}, this cardinality is bounded below by $\left(\sum_{i=1}^h |C(p,q_i)|\right)-2h+3$, which completes the proof.
\end{proof}

\begin{theorem} Let $h\geq 1$ and suppose $p,q_1,\ldots,q_h$ are distinct vertices of $\spts$. Then, there exists a family of sets $J$ such that:
	\begin{itemize}
		\item For all $X\in J$, $X$ is a nonempty subset of $\spts$ in convex position.
		\item If $\{X, X'\}\subseteq J$ and $X\neq X'$, then $\mbox{\emph{conv}}(X)\cap\mbox{\emph{conv}}(X')=\emptyset$.
		\item $|J|\leq 4h-2$.
		\item $\sum_{X\in J} |X|\geq \left(\frac{1}{2}\sum_{i=1}^h |C(p,q_i)|\right)-2h+2$.
	\end{itemize} 
	\label{convex2_theorem}
\end{theorem}

\begin{proof} If $h=1$, then $J=\{C(p,q_1)\}$ satisfies the conclusion of the theorem. Now suppose $h>1$. Without loss of generality, assume $q_1, \ldots, q_h$ are listed in counterclockwise order about $p$. If $q_i\in\mc{H}^+(p,q_1)$ for all $2\leq i\leq h$, the result follows immediately from Corollary~\ref{convex1_corollary}. If not, let $j$ be the minimum $i$ such that $q_{i+1}\in\mc{H}^-(p,q_1)$. Clearly, $1\leq j<h$, and $q_i\in\mc{H}^+(p,q_1)$ for $2\leq i\leq j$, and $q_i\in\mc{H}^-(p,q_1)$ for $j+1\leq i\leq h$. 
	
	First, let us suppose $\sum_{i=1}^j|C(p,q_i)|\geq \frac{1}{2}\sum_{i=1}^h|C(p,q_i)|$. If $j=1$, then $J=\{C(p,q_1)\}$ satisfies the conclusions of the theorem. If $j>1$, we apply Corollary~\ref{convex1_corollary} to the points $p,q_1,\ldots,q_j$ to obtain a family $J$ of nonempty subsets in convex position with disjoint convex hulls such that $|J|\leq 4j-2\leq 4h-2$ and
	$\sum_{X\in J} |X|\geq \left(\sum_{i=1}^j |C(p,q_i)|\right)-2j+3\geq \left(\frac{1}{2}\sum_{i=1}^h|C(p,q_i)|\right)-2h+2$.
	
	Now suppose $\sum_{i=1}^j|C(p,q_i)|<\frac{1}{2}\sum_{i=1}^h|C(p,q_i)|$. This implies $\sum_{i=j+1}^h|C(p,q_i)|>\frac{1}{2}\sum_{i=1}^h|C(p,q_i)|$. If $j+1=h$, then $J=\{C(p,q_{j+1})\}$ satisfies the conclusions of the theorem. If $j+1<h$, note that $q_i\in\mc{H}^+(p,q_{j+1})$ for $j+2\leq i\leq h$ because $q_i\in\mc{H}^-(p,q_1)$ for $j+1\leq i\leq h$. Applying Corollary~\ref{convex1_corollary} to the points $p,q_{j+1},\ldots,q_h$, we obtain a family $J$ of nonempty subsets in convex position with disjoint convex hulls such that $|J|\leq 4(h-j)-2\leq 4h-2$ and
	$\sum_{X\in J} |X|\geq \left(\sum_{i=j+1}^h |C(p,q_i)|\right)-2(h-j)+3\geq \left(\frac{1}{2}\sum_{i=1}^h|C(p,q_i)|\right)-2h+2$.
\end{proof}

\section{Empty Pentagons and Hexagons in Convex Subsets}
\begin{lemma} Let $m$ and $k$ be integers larger than $1$, and let $D^{\circ}$ be any set. Suppose $w:D^{\circ}\rightarrow\mathbb{R}$ satisfies the following conditions: 
	\begin{itemize}
		\item $w(x)\in\left\{\frac{1}{k-1},\frac{1}{k-2},\ldots,\frac{1}{2},1\right\}$, for all $x\in D^{\circ}$.
		\item $\left|w^{-1}(1/i)\right|\leq 2m$, for all $1\leq i\leq k-1$.
	\end{itemize}
	Then, $\sum_{x\in D^{\circ}} w(x)> 2m\ln\!\left(\frac{2mk}{2m(k+1)-|D^{\circ}|}\right)$.
	\label{weight_function_lemma}
\end{lemma}

\begin{proof} Note that the conclusion is trivially true for $|D^{\circ}|<2m$ since the right side of the inequality is negative. For brevity, let $t=|D^{\circ}|$. Assume $t\geq 2m$.  Clearly, $t=\sum_{i=1}^{k-1}|w^{-1}(1/i)|\leq\sum_{i=1}^{k-1} 2m=2m(k-1)$, by the conditions on $w$. 
	
	Consider the sequence $s=\left(\frac{1}{k-1},\ldots,\frac{1}{k-1},\frac{1}{k-2},\ldots,\frac{1}{k-2},\ldots,1,\ldots,1\right)$ in which each of the values $1/i$ is repeated $2m$ times. Note that the conditions on $w$ imply that the terms of $\sum_{x\in D^{\circ}}w(x)$, when sorted in non-decreasing order, form a subsequence of $s$ with length $t$. Clearly, the sum of the terms in this subsequence is greater than or equal to the sum of the first $t$ terms of $s$ (i.e., the $t$ smallest terms of $s$): 
	$\sum_{x\in D^{\circ}}w(x)\geq \sum_{i=1}^{t}s_{i}.$
	 
	On the other hand, 
	\begin{align*}
		\sum_{i=1}^{t}s_{i}\geq \sum_{i=1}^{2m\lfloor t/(2m)\rfloor}s_i&=2m\sum_{j=1}^{\lfloor t/(2m)\rfloor}\tfrac{1}{k-j}\\
		&=2m\sum_{j=k-\lfloor t/(2m)\rfloor}^{k-1} 1/j\geq 2m\left(\ln(k)-\ln(k-\lfloor t/(2m)\rfloor)\right)
	\end{align*}
	since $\sum_{j=a}^b\frac{1}{j}\geq\ln(b+1)-\ln(a)$. But $k-\lfloor t/(2m)\rfloor<k+1-t/(2m)$. Therefore:
	\[\sum_{x\in D^{\circ}}w(x)
	> 2m\left(\ln(k)-\ln(k+1-t/(2m))\right)
	= 2m\ln\!\left(\tfrac{2mk}{2m(k+1)-t}\right).\]
\end{proof}

\begin{theorem} 
	Suppose $B\subseteq \spts$ is in convex position and is nonempty. Let $m=|B|$ and let $n_0$ be the number of points of $\spts$ in the interior of $\conv(B)$, i.e., $n_0=|(\spts\cap\conv(B))-B|$. Finally, let $\ps\!_0=\{Q\in \ps: Q\subseteq \conv(B)\}$. Then,  
	\[|\ps\!_0|>\tfrac{1}{5}m^2e^{-n_0/m}-\tfrac{11}{5}m.\] 
	Additionally, for every $\varepsilon>0$, there exists $\delta>0$ such that if $n_0/m<\ln(2)-\varepsilon$, then the number of empty convex hexagons in $\conv(B)$ is at least $\delta m^2$, for sufficiently large $m$.
	
	\label{pentagons_in_convex_set_Theorem}
\end{theorem}

\begin{proof} The result is trivially true for $1\leq m<11$ since $\tfrac{1}{5}m^2e^{-n_0/m}-\tfrac{11}{5}m \leq \tfrac{1}{5}m(m-11) <0$ in that case. Now suppose $m\geq 11$ and let $B_0=(\spts\cap\conv(B))-B$, the points of $\spts$ in the interior of $\conv(B)$, so that $n_0=|B_0|$. 
	
	Starting at any point $z_0$ of $B$, let $z_0,\ldots,z_{m-1}$ be the list of all points in $B$ as they appear when we traverse the boundary of $\conv(B)$ in counterclockwise order. Let $E_0=\{(z_i,z_{i+j}): 0\leq i\leq m-1 \text{ and } 3\leq j\leq m-3\}$, with the understanding that the operation $i+j$ in the subscript is calculated modulo $m$. Note that for all $(p,q)\in E_0$, $\left|B\cap\mc{H}^+(p,q)\right|\geq 2$ and $\left|B\cap\mc{H}^-(p,q)\right|\geq 2$. Also, $|E_0|=m(m-5)$. 
	
	We say that $(z,z')\in E_0$ is blocked (in $B$) if its visibility set contains only one point in $B\cup B_0$, i.e., $\left|V(z,z')\cap(B\cup B_0)\right|=1$. If $(z,z')$ is blocked and $V(z,z')\cap(B\cup B_0)=\{q\}$, we say that $(z,z')$ is blocked by $q$. Clearly, $(z,z')$ is blocked by $q$ if and only if $q\in\conv(z,z',w)$ for all $w\in(B\cup B_0)\cap\mc{H}^+(z,z')$. Note that if $(z,z')$ is blocked by $q$, then necessarily $q\in B_0$ because $B$ is in convex position. See Figure~\ref{Figure13} for an example. 
	
	\begin{figure}[tpb!]
		\centering\includegraphics[scale=.75]{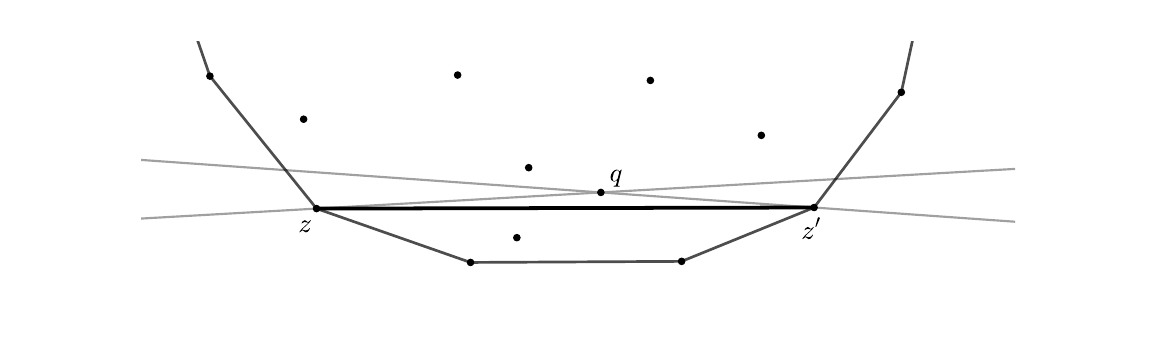}   
		\caption{$(z,z')$ is blocked by $q$.}
		\label{Figure13}
	\end{figure}
	
	Let $D^{\circ}$ be the set of edges in $E_0$ that are blocked. Now we show that if an edge $(z,z')$ is \emph{not} blocked, then there exists a pentagon in $\ps\!_0$ with $(z,z')$ as a diagonal. Let $(z,z')\in E_0-D^\circ$. Then, there exist $q_1$ and $q_2$ in $V(z,z')\cap(B\cup B_0)$ such that $q_1\neq q_2$. Clearly, $\conv(z,z',q_i)$ is empty (for $i=1,2$) and $\conv(z,z',q_1,q_2)$ is a convex quadrilateral. Now select $q_2'$ in $\conv(z,z',q_1,q_2)$ such that $\conv(z,z',q_1,q_2')$ is an empty convex quadrilateral. Then we select $q_3$ in $(B\cup B_0)\cap\mc{H}^-(z,z')$ such that $\conv(z,q_3,z')$ is empty. Finally, define $Q(z,z')$ to be $\conv(z,q_3,z',q_1,q_2')$. Since $Q(z,z')=\conv(z,q_3,z')\cup\conv(z,z',q_1,q_2')$, $Q(z,z')$ is empty. Moreover, $Q(z,z')$ is convex because $z$ and $z'$ lie on the boundary of $\conv(B)$ and $Q(z,z')\subseteq\conv(B)$. Therefore, $Q(z,z')$ is an empty convex pentagon, i.e., $Q(z,z')\in \ps\!_0$. The function $(z,z')\rightarrow Q(z,z')$ is at most five-to-one because $(z,z')$ is a diagonal of $Q(z,z')$ with two points on its positive side, and there are five such diagonals in $Q(z,z')$. Therefore, 
	\[5|\ps\!_0|\geq|E_0-D^\circ|=|E_0|-|D^\circ|=m(m-5)-|D^\circ|.\]
	
	Additionally, note that if both $(z,z')$ and $(z',z)$ are not blocked, then the procedure above can be adapted to produce an empty convex hexagon with two opposite vertices in $B$. Since there are three pairs of opposite vertices in any hexagon, it follows that the number of empty convex hexagons in $\conv(B)$ is at least one third the number of 2-sets $\{z,z'\}\subset B$ such that both $(z,z')$ and $(z',z)$ are not blocked.
	
	To bound $|D^\circ|$, we consider the properties of edges in $D^\circ$ that are blocked by the same point. First, we show that two such edges must intersect, i.e., they must cross or have a common endpoint. Suppose $(z,z')$ and $(y,y')$ are both blocked by $q$ and they do not intersect. By the convexity of $B$, the points $z,z',y,y'$ must appear in one of the following four ways along the boundary of $\conv(B)$ when listed in counterclockwise order starting at $z$: $z,z',y,y'$ or $z,z',y',y$ or $z,y,y',z'$ or $z,y',y,z'$ (since any order alternating $z$'s and $y$'s would make the edges cross). In the first case (order $z,z',y,y'$), we have $\{y,y'\}\subset\mc{H}^+(z,z')$ and $\{z,z'\}\subset\mc{H}^+(y,y')$. Therefore, $q\in\conv(z,z',y)$ since $q$ blocks $(z,z')$, which implies $q\in\mc{H}^-(z,y)$. On the other hand, $q\in\conv(y,y',z)$ since $q$ blocks $(y,y')$, which implies $q\in\mc{H}^-(y,z)=\mc{H}^+(z,y)$, so this case is contradictory. In the second case (order $z, z', y', y$), we have $\{y,y'\}\subset\mc{H}^+(z,z')$ and $\{z,z'\}\subset\mc{H}^-(y,y')$. Hence, $q\in\conv(z,z',y')$ since $q$ blocks $(z,z')$, which implies $q\in\mc{H}^-(y,y')$. But $q$ has to be in $\mc{H}^+(y,y')$ to block $(y,y')$, so this case is contradictory too. The third and fourth cases are similar. Therefore, two edges that are blocked by the same point must intersect.
	
	\begin{figure}[tpb!]
		\centering\includegraphics[scale=.475]{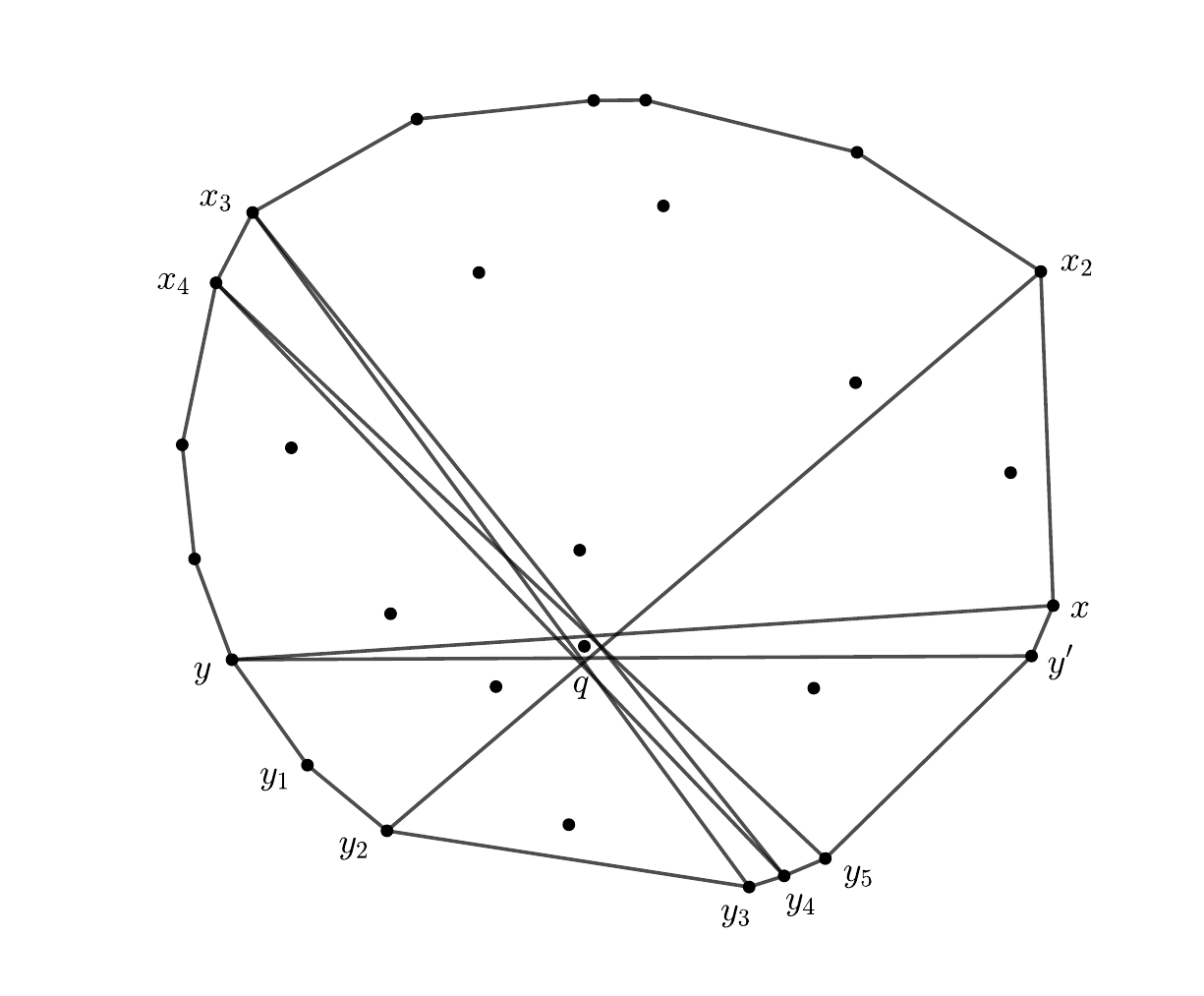}  
		\caption{$(x,y),(y,y'),(y_2,x_2),(y_5,x_4),(x_4,y_4),(y_4,x_3)$ and $(x_3,y_3)$ are blocked by $q$.}
		\label{Figure14}
	\end{figure}
		
	However, it is clear that two edges with the same left endpoint, say $(z,z')$ and $(z,z'')$, cannot be blocked by the same point $q$. Indeed, assume without loss of generality that $z''\in\mc{H}^+(z,z')$ and $q$ blocks $(z,z')$. Then, necessarily $q\in\conv(z,z',z'')$, which implies $q\in\mc{H}^-(z,z'')$, so $q$ cannot block $(z,z'')$. The same reasoning applies to edges with the same right endpoint. Therefore, any two edges with a common endpoint $z'$ that are blocked by the same point $q$ must be of the form $(z,z')$ and $(z',z'')$. It follows that, for every $z'\in B$ and $q\in B_0$, the maximum number of edges incident to $z'$ that are blocked by $q$ is two.
	
	Now let $f:D^\circ\rightarrow B_0$ be defined by $f(z,z')=q$, where $q$ is the unique point in $B_0$ that blocks $(z,z')$. Also, define the length $l(z,z')$ to be the minimum of $|B\cap\mc{H}^-(z,z')|$ and $|B\cap\mc{H}^+(z,z')|$ for any $(z,z')\in E_0$. Finally, define the weight function $w(z,z')$ to be $1/(2+l(z,z'))$. Note that $2\leq l(z,z')\leq \lfloor \frac{m}{2}\rfloor-1$ and $1/(\lfloor\frac{m}{2}\rfloor+1)\leq w(z,z')\leq 1/4$ for all $(z,z')\in E_0$.
	
	We now show that $|f^{-1}(q)|\leq 3+2\min\{l(z,z'): (z,z')\in f^{-1}(q)\}$ for any $q\in f(D^\circ)$. Indeed, let $(y,y')\in f^{-1}(q)$ be such that $l(y,y')=\min\{l(z,z'): (z,z')\in f^{-1}(q)\}$. As explained above, every edge in $f^{-1}(q)$ must either cross $(y,y')$ or be incident to $y$ or to $y'$. Furthermore, other than $(y,y')$, there may exist at most one edge in $f^{-1}(q)$ that is incident to $y$ and one that is incident to $y'$. Thus, the total number of edges in $f^{-1}(q)$ incident to either $y$ or $y'$ is at most three. All other edges in $f^{-1}(q)$ must cross $(y,y')$, i.e., all other edges must have an endpoint in $B\cap\mc{H}^-(y,y')$ and one in $B\cap\mc{H}^+(y,y')$. But either $B\cap\mc{H}^-(y,y')$ or $B\cap\mc{H}^+(y,y')$ has exactly $l(y,y')$ vertices, and each of them can be incident to at most two edges in $f^{-1}(q)$. Therefore, 
	\[\left|f^{-1}(q)\right|\leq 3+2l(y,y')=3+2\min\{l(z,z'): (z,z')\in f^{-1}(q)\}.\] 
	See Figure~\ref{Figure14} for an example.
		
	Let $q\in f(D^\circ)$. Denote $\max\{w(z,z'):(z,z')\in f^{-1}(q)\}$ by $w_{max}$ and $\min\{l(z,z'):(z,z')\in f^{-1}(q)\}$ by $l_{min}$. Clearly, $w_{max}=1/(2+l_{min})$. Therefore, 
	\[\sum_{(z,z')\in f^{-1}(q)}w(z,z')\leq\sum_{(z,z')\in f^{-1}(q)}w_{max}=|f^{-1}(q)|w_{max}=|f^{-1}(q)|\left(\tfrac{1}{2+l_{min}}\right)\leq\tfrac{3+2l_{min}}{2+l_{min}}<2.\]
	
	Hence, $\sum_{(z,z')\in f^{-1}(q)}w(z,z')<2$ for all $q\in f(D^\circ)$. Summing over all points $q$ yields:
	\[\sum_q\sum_{(z,z')\in f^{-1}(q)}w(z,z')<2|f(D^\circ)|\leq 2|B_0|=2n_0.\] 
	Because $\sum_q\sum_{(z,z')\in f^{-1}(q)}w(z,z')=\sum_{(z,z')\in D^\circ}w(z,z')$, it follows that $\sum_{(z,z')\in D^\circ}w(z,z')<2n_0$, i.e., the sum of the weights of all the blocked edges is at most $2n_0$.
	
	On the other hand, $m$, $D^\circ$, and $w$ satisfy the conditions of Lemma~\ref{weight_function_lemma} for $k=\lfloor \frac{m}{2}\rfloor +2$. Therefore:
	\begin{align*}
		\sum_{(z,z')\in D^\circ} w(z,z') &> 2m\ln\!\left(\tfrac{2m(\lfloor \frac{m}{2}\rfloor +2)}{2m\left((\lfloor \frac{m}{2}\rfloor +2)+1\right)-|D^{\circ}|}\right)\\
		&> 2m\ln\!\left(\tfrac{m(2\lfloor \frac{m}{2}\rfloor +4)}{m\left(2\lfloor \frac{m}{2}\rfloor +6\right)-|D^{\circ}|}\right)\\
		&> 2m\ln\!\left(\tfrac{m^2}{m(m+6)-|D^{\circ}|}\right).
	\end{align*}
	
	Combining the two inequalities yields $2n_0>2m\ln\!\left(\frac{m^2}{m(m+6)-|D^{\circ}|}\right)$. Solving for $|D^\circ|$ we obtain a bound on the maximum number of edges from $E_0$ that can be blocked with $n_0$ points: 
	\[|D^{\circ}|<m^2(1-e^{-n_0/m})+6m.\]
	
	But we established above that $5|\ps\!_0|\geq m(m-5)-|D^\circ|$. Substituting, we obtain the desired bound on pentagons:
	\[5|\ps\!_0| > m(m-5)-m^2(1-e^{-n_0/m})-6m = m^2e^{-n_0/m}-11m,\] 
	so $|\ps\!_0|>\frac{1}{5}m^2e^{-n_0/m}-\frac{11}{5}m$.
	
	Finally, let $\varepsilon>0$ and suppose $n_0/m<\ln(2)-\varepsilon$. Then $1-e^{-n_0/m}<1-e^{\varepsilon}/2$. Applying this to the bound on $|D^{\circ}|$ above yields: 
	\[|D^{\circ}|< m^2(1-e^{\varepsilon}/2)+6m,\] 
	so $m(m-5)/2-|D^{\circ}|>\frac{1}{2}(e^{\varepsilon}-1)m^2-17m/2>\frac{1}{4}(e^{\varepsilon}-1)m^2$, for sufficiently large $m$.
	
	But the number of pairs $\{z,z'\}$ such that both $(z,z')$ and $(z',z)$ are in $E_0-D^{\circ}$ is at least $m(m-5)/2-|D^{\circ}|$. Letting $\delta=\frac{1}{12}(e^{\varepsilon}-1)$, we obtain that if $n_0/m<\ln(2)-\varepsilon$, then the number of empty convex hexagons in $\conv(B)$ is at least $\delta m^2$ for sufficiently large $m$. 
\end{proof}

\begin{corollary} 
	Let $k\geq 1$ and suppose $B_i\subseteq \spts$ is a nonempty subset of $\spts$ in convex position for $1\leq i\leq k$. Additionally, suppose $\conv(B_i)\cap \conv(B_j)=\emptyset$ for $i\neq j$. Let $m_i=|B_i|$, $n_i=|(\spts\cap\conv(B_i))-B_i|$, $m=\sum_{i=1}^k m_i$, and $n_0=\sum_{i=1}^k n_i$. Then, 
	\[|\ps|>\tfrac{1}{5k}m^2e^{-n_0/m}-\tfrac{11}{5}m.\] 
	\label{pentagons_in_convex_sets_Corollary}
\end{corollary}

\begin{proof} Let $f:\mathbb{R}^+\times\mathbb{R}\rightarrow \mathbb{R}$ be defined by $f(x,y)=x^2e^{-y/x}-11x$. Since 
	$\frac{\partial^2 f}{\partial y^2}=e^{-y/x}>0 \quad \text{and} \quad \frac{\partial^2 f}{\partial x^2}\frac{\partial^2 f}{\partial y^2}-\left(\frac{\partial^2 f}{\partial x\partial y}\right)^2=e^{-2y/x}>0$, 
	the function $f$ is convex. Therefore, by Jensen's inequality:
	\[\tfrac{1}{k}\sum_{i=1}^k f(m_i,n_i)\geq f\!\left(\tfrac{1}{k}\sum_{i=1}^k m_i,\,\tfrac{1}{k}\sum_{i=1}^k n_i\right)=f\left(\tfrac{m}{k},\tfrac{n_0}{k}\right)=\left(\tfrac{m}{k}\right)^2e^{-\tfrac{n_0}{k}/(\tfrac{m}{k})}-11\tfrac{m}{k}\] 
	Multiplying by $k$ and simplifying, we obtain: 
	$\sum_{i=1}^k f(m_i,n_i)\geq\frac{1}{k}m^2e^{-n_0/m}-11m$
	 
	On the other hand, by Theorem~\ref{pentagons_in_convex_set_Theorem}, $5|\ps\!_i|>f(m_i,n_i)$, where $\ps\!_i=\{Q\in \ps: Q\subseteq\conv(B_i)\}$ for $1\leq i\leq k$. Summing this over all $k$ subsets gives $5\sum_{i=1}^k|\ps\!_i|>\sum_{i=1}^k f(m_i,n_i)$. Finally, $5|\ps|\geq \sum_{i=1}^k 5|\ps\!_i|$, since $\conv(B_i)\cap\conv(B_j)=\emptyset$ for $i\neq j$. Therefore, 
	$5|\ps|>\tfrac{1}{k}m^2e^{-n_0/m}-11m$, so $|\ps|>\tfrac{1}{5k}m^2e^{-n_0/m}-\tfrac{11}{5}m$.
\end{proof}

\section{The $\Omega(n^{3/2})$ Lower Bound for Pentagons}
We now have all the necessary components for the general case. Combining  the sets in convex position from Theorem~\ref{convex2_theorem}, the bound on pentagons for the convex case from Corollary~\ref{pentagons_in_convex_sets_Corollary}, and the bounds from Theorem~\ref{bound_on_edges_theorem}, we can finally establish the $\Omega(n^{3/2})$ lower bound on the number of empty pentagons for arbitrary sets in general position with sufficiently many points.

\begin{theorem} 
	There exist constants $K>0$, $N_0>0$ such that $|\ps|> Kn^{3/2}$ for $n>N_0$.
	\label{Pentagons_Theorem}
\end{theorem}

\begin{proof} Let $K=\frac{1}{2^{800}}$ and $N_0=2^{1603}$. Suppose $n>N_0$. 
	
	By Theorem~\ref{bound_on_edges_theorem}, there exists $D\subseteq \es$ satisfying the following two properties:
	\begin{itemize}
		\item $|\ps|\geq \tfrac{1}{20}|D|$ 
		\item $6|\ps|+\sum_{(p,q)\in D}|C(p,q)|\geq \tfrac{1}{16}n^2-\tfrac{17}{16}n$
	\end{itemize}
	We consider two cases. First, suppose that either $|D|\geq n^{3/2}$ or $|\ps|\geq\frac{1}{192}n^2$. Then, by the first property of $D$ above, either $|\ps|\geq \frac{1}{20}n^{3/2}$ or $|\ps|\geq\frac{1}{192}n^2$. Therefore, $|\ps|\geq Kn^{3/2}$ and we are done. 
	
	For the second case, suppose $|D|< n^{3/2}$ and $|\ps|<\frac{1}{192}n^2$. Then, by the second property of $D$ above, 
	$\sum_{(p,q)\in D}|C(p,q)|\geq\frac{1}{16}n^2-\frac{17}{16}n-6|\ps|\geq \frac{1}{16}n^2-\frac{17}{16}n-\frac{6}{192}n^2=\frac{1}{32}n^2-\frac{17}{16}n$. 
	Therefore, we have: 
	\[\sum_{(p,q)\in D}|C(p,q)|\geq\tfrac{1}{32}n^2-\tfrac{17}{16}n.\] 
	
	Now we group the edges in $D$ according to their left endpoints. Let $D_p=\{q\in \spts: (p,q)\in D\}$. Clearly, there is a one-to-one correspondence between the edges of $D$ left-incident to $p$ and the vertices in $D_p$. Hence, $\sum_{p\in \spts}|D_p|=|D|<n^{3/2}$. We say that a vertex $p$ has high degree if $|D_p|>256\sqrt{n}$.  
	Let $\spts_h$ be the set of vertices of $\spts$ with high degree, i.e., $\spts_h=\{p\in \spts: |D_p|>256\sqrt{n}\}$. 
	Then, $\sum_{p\in \spts_h}|D_p|>\sum_{p\in \spts_h}256\sqrt{n}=256|\spts_h|\sqrt{n}$. But $n^{3/2}>\sum_{p\in \spts}|D_p|\geq\sum_{p\in \spts_h}|D_p|$. Therefore, $n^{3/2}> 256|\spts_h|\sqrt{n}$, or $|\spts_h|<n/256$. 
	
	On the other hand, by Corollary~\ref{upperboundCCC_corollary}, $\sum_{q\in D_p}|C(p,q)|\leq 4(n-1)< 4n$ for any $p\in \spts$. Therefore, 
	$\sum_{p\in \spts_h}\sum_{q\in D_p}|C(p,q)|< \sum_{p\in \spts_h} 4n= 4n|\spts_h|< 4n\left(\frac{n}{256}\right)=\frac{1}{64}n^2.$ 
	Hence, we have: 
	\[\sum_{p\in \spts_h}\sum_{q\in D_p}|C(p,q)|<\tfrac{1}{64}n^2.\]
	
	Let $\spts_l=\{p\in \spts: |D_p|\leq 256\sqrt{n}\}$, the low-degree vertices of $\spts$. Since $\spts_h$ and $\spts_l$ partition $\spts$ and $\sum_{p\in \spts_h}\sum_{q\in D_p}|C(p,q)|<\frac{1}{64}n^2$, it follows that 
	$\frac{1}{64}n^2+\sum_{p\in \spts_l}\sum_{q\in D_p}|C(p,q)|> \sum_{(p,q)\in D}|C(p,q)|\geq\frac{1}{32}n^2-\frac{17}{16}n$. 
	Therefore, $\sum_{p\in \spts_l}\sum_{q\in D_p}|C(p,q)|> \tfrac{1}{32}n^2-\tfrac{17}{16}n-\tfrac{1}{64}n^2$, so we have: 
	\[\sum_{p\in \spts_l}\sum_{q\in D_p}|C(p,q)|> \tfrac{1}{64}n^2-\tfrac{17}{16}n.\]
	Using the obvious bound $|\spts_l|\leq n$, we obtain that the average value of $\sum_{q\in D_p}|C(p,q)|$ over $p\in \spts_l$ satisfies:
	\[\frac{1}{|\spts_l|}\sum_{p\in \spts_l}\sum_{q\in D_p}|C(p,q)|>  \tfrac{1}{n}\left(\tfrac{1}{64}n^2-\tfrac{17}{16}n\right) > \tfrac{1}{128}n\] 
	since $n>N_0$. Therefore, there exists an element $\hat{p}\in \spts_l$ such that $\sum_{q\in D_{\hat{p}}}|C(\hat{p},q)|>\frac{1}{128}n$.
	
	Let $h=\left|D_{\hat{p}}\right|$ and $D_{\hat{p}}=\{q_1,\ldots,q_h\}$. By Theorem~\ref{convex2_theorem}, there exists a family of sets $\{B_1,\ldots,B_k\}$ such that $B_i$ is a nonempty subset of $\spts$ in convex position for all $i$, $\conv(B_i)\cap\conv(B_j)=\emptyset$ for all $i\neq j$, and 
	$\sum_{i=1}^k |B_i|\geq \left(\frac{1}{2}\sum_{i=1}^h|C(\hat{p},q_i)|\right)-2h+2$, where $1\leq k\leq 4h-2$. Note that $\sum_{i=1}^h|C(\hat{p},q_i)|>\frac{1}{128}n$, by definition of $\hat{p}$, and $h=\left|D_{\hat{p}}\right|\leq 256\sqrt{n}$, since $\hat{p}\in \spts_l$. Therefore, $k<1024\sqrt{n}$ and    
	\[\sum_{i=1}^k|B_i|> \tfrac{1}{2}\!\left(\tfrac{1}{128}n\right)-2(256\sqrt{n})+2>\tfrac{1}{256}n-512\sqrt{n}>\tfrac{1}{512}n\]
	since $n>N_0$.
	
	Let $m_i=|B_i|$, $n_i=|(\spts\cap\conv(B_i))-B_i|$, $m=\sum_{i=1}^k m_i$, and $n_0=\sum_{i=1}^k n_i$. By Corollary~\ref{pentagons_in_convex_sets_Corollary}, $|\ps|>\frac{1}{5k}m^2e^{-n_0/m}-\tfrac{11}{5}m$. We have shown that $m>\frac{1}{512}n$ and $k<1024\sqrt{n}$. Also, clearly $n_0\leq n$ and $m\leq n$. Therefore, 
	\begin{align*}
		|\ps| &> \tfrac{1}{5k}m^2e^{-n_0/m}-\tfrac{11}{5}m\\
		&> \tfrac{1}{5(1024\sqrt{n})}\left(\tfrac{1}{512}n\right)^2e^{-n/(n/512)}-\tfrac{11}{5}n\\
		&> \tfrac{e^{-512}}{8(2^{28})}n^{3/2}-\tfrac{11}{5}n\\
		&> \tfrac{1}{2^{799}}n^{3/2}-2\sqrt{2}n\\
		&> \tfrac{1}{2^{800}}n^{3/2}
	\end{align*}
	since $e<2^{3/2}$ and $n>N_0=2^{1603}$.
\end{proof}

Clearly, no attempt was made to optimize the values of the constants $K$ and $N_0$ in the preceding proof. While this $\Omega(n^{3/2})$ bound represents a significant improvement over previously published bounds, the true order is widely conjectured to be $\Theta(n^2)$, so we consider that there is little value in improving these constants.


\end{document}